\documentclass[a4paper,10pt]{article}

\usepackage[margin=0.7in]{geometry}
\usepackage[english]{babel}
\usepackage{amsfonts}
\usepackage{amssymb}
\usepackage{amsthm, amsmath}
\usepackage{booktabs}
\usepackage{eucal}
\usepackage{mathrsfs}
\allowdisplaybreaks[0]
\usepackage{bbm}
\usepackage{graphicx}
\usepackage{caption}
\usepackage{subcaption}
\usepackage{algorithm,algorithmic}
\usepackage{framed}
\usepackage{ftnxtra}
\usepackage[titletoc,title]{appendix}
\usepackage{booktabs}
\usepackage[hyphens]{url}
\usepackage[hypertexnames=false]{hyperref}

\usepackage[usenames,dvipsnames]{color}

\def\red #1{{\color{red}#1}}

\let\red\relax

\newcommand{\vz}[1]{\ensuremath{\mathbb{#1}}}
\newcommand{\Z}{{\vz Z}}
\newcommand{\N}{{\vz N}}
\newcommand{\R}{{\vz R}}

\newtheorem{theorem}{Theorem}[section]

\newtheorem{remark}[theorem]{Remark}

\numberwithin{equation}{section}

\providecommand{\norm}[1]{\left\lVert#1\right\rVert}

\numberwithin{equation}{section}

\begin{document}

\thispagestyle{empty}

\begin{center}
{\Large \bf 
Graph clustering, variational image segmentation methods and \\Hough transform scale detection for object measurement in images} \\
\end{center}

\begin{center}
     \footnotesize{            {\sc Luca Calatroni}\\
 MIDA group, Dipartimento di Matematica, Universit\'a degli studi di Genova, \\ Via Dodecaneso 35, 16146, Italy,\\
 (calatroni@dima.unige.it) \\[0.4cm]
 {\sc Yves van Gennip}\\
 School of Mathematical Sciences, The University of Nottingham\\
 University Park, NG7 2RD, Nottingham, UK \\
 (y.vangennip@nottingham.ac.uk) \\[0.4cm]
  {\sc Carola-Bibiane Sch\"{o}nlieb} \\
 Department of Applied Mathematics and Theoretical Physics (DAMTP)\\
 University of Cambridge, Wilberforce Road, CB3 0WA, Cambridge, UK \\
 (cbs31@cam.ac.uk)} \\[0.4cm]
  {\sc Hannah Rowland}\\
 Department of Zoology, University of Cambridge, Downing Street,\\ CB2 3EJ, Cambridge, UK  \\
  Institute of Zoology, Zoological Society of London, Regents Park, \\
   NW1 4RY, London, UK \\
(hr325@cam.ac.uk)\\[0.4cm]
  {\sc Arjuna Flenner}\\
Image and Signal Processing Branch, NAVAIR\\
1 Administrative Circle, China Lake CA, USA\\
(arjuna.flenner@navy.mil)
       \end{center}

\begin{abstract}
We consider the problem of scale detection in images where a region of interest is present together with a measurement tool (e.g. a ruler). For the segmentation part, we focus on the graph based method presented in \cite{BerFlen} which reinterprets classical continuous Ginzburg-Landau minimisation models in a totally discrete framework. To overcome the numerical difficulties due to the large size of the images considered we use matrix completion and splitting techniques. The scale on the measurement tool is detected via a Hough transform based algorithm. The method is then applied to some measurement tasks arising in real-world applications such as zoology, medicine and archaeology. 
\end{abstract}

\section{Introduction} \label{sec:int}

\emph{Image segmentation} denotes the task of partitioning an image in its constituent parts. Feature based segmentation looks at distinctive characteristics (features) in the image, grouping similar pixels into clusters which are meaningful for the application at hand. Typical examples of features are based on greyscale/RGB intensity and texture. Mathematical methods for image segmentation are mainly formalised in terms of variational problems in which the segmented image is a minimiser of an energy. The most common image feature encoded in such energies is the magnitude of the image gradient, detecting regions (or contours) where sharp variations of the intensity values occur. 
Examples include the Mumford-Shah segmentation approach \cite{mumfshah}, the snakes and geodesic active contour models \cite{KaasSnakes,CasellesGeodesic}. Moreover, in \cite{chanvese} Chan and Vese proposed an instance of the Mumford-Shah model for piecewise constant  images whose energy is based on the mean greyvalues of the image inside and outside of the segmented region rather than the image gradient and hence does not require strong edges for segmentation. \red{The Chan-Vese model has been extended for vector-valued images such as RGB images in \cite{chansandvese}.} Other image segmentation methods have been considered in \cite{kohn,esedoglu}.  They rely on the use of the total variation (TV) seminorm \cite{ambrosio}, which is commonly used for image processing tasks due to its properties of simultaneous edge preservation and smoothing (see \cite{ROF}).

The non-smoothness of most of the segmentation energies renders their numerical minimisation usually difficult. In the case of the Mumford-Shah segmentation model the numerical realisation is additionally complicated by its dependency on the image function as well as the object contour. To overcome this, several regularisation methods and approximations have been proposed in the literature, e.g. \cite{ambrosiotort,braides,braidesdalmaso,chanvese1} for Mumford-Shah segmentation. In the context of TV based segmentation models the Ginzburg-Landau functional has an important role. Originally considered for the modelling of physical phenomena such as phase transition and phase separation (cf. \cite{brokate} for a survey on the topics) it is used in imaging for approximating the TV energy. Some examples of the use of this functional in the context of image processing are \cite{esedoglu,esed1,esed2}, which relate to previous works by Ambrosio and Tortorelli on diffuse interface approximation models \cite{ambrosio,ambrosiotort}. 



Such variational methods for image segmentation have been extensively studied from an analytical point of view and the segmentation is usually robust and computationally efficient. However, variational image segmentation as described above still faces many problems in the presence of low contrast and the absence of clear boundaries separating regions. Their main drawback is that they are limited to image features which can be mathematically formalised (e.g. in terms of an image gradient) and encoded within a segmentation energy. In recent years dictionary based methods have become more and more popular in the image processing community, complementing more classical variational segmentation methods. By learning the distinctive features of the region to be segmented from examples provided by the user, these methods are able to segment the desired regions in the image correctly.  

\medskip

In this work, we consider the method proposed in \cite{BerFlen,BerFlen1,BerFlen2,MerkujevBert2015} for image segmentation and labelling. This approach goes beyond the standard variational approach in two respects. Firstly, the model is set up in the purely discrete framework of graphs. This is rather unusual for variational models where one normally considers functionals and function spaces defined on subdomains of $\R^2$ in order to exploit properties and tools from convex and functional analysis and calculus of variations. 
Secondly, the new framework allows for more flexibility in terms of the features considered. Additional features like texture, light intensity or others, can be considered as well without encoding them in the function space or the regularity of the functions. Due to the possibly very large size of the image (nowadays of the order of megapixel for professional cameras) and the large number of features considered, the construction of the problem may be computationally expensive and often requires reduction techniques \cite{nystrom,naeini,fowlkes}.
In several papers (see, e.g., \cite{shimalik,stoer,guattery}) the segmentation problem was rephrased in the graph framework by means of the graph cut objective function. \red{Follow-up works on the use of graph-based approaches are, for instance, \cite{MerkujevKostic2013,MerkujevSunu2014} where an iterative application of heat diffusion and thresholding, also known as the  Merriman-Bence-Osher (MBO) method \cite{MBO1992} is discussed for binary image labelling, and \cite{hu} where the Mumford-Shah model is reinterpreted in a graph setting.}

In this paper, we also address the problem of detection of objects with geometrical properties that are \emph{a priori} known. An example is the detection of lines and circles. These objects can be identified by mapping them onto an auxiliary space where relevant geometrical properties (such as linear alignment and roundness) are represented as peaks of specific auxiliary functions. In this work, we use the Hough transform \cite{hough} to detect measurement tools (rulers, concentric circles of fixed radii) with the intent of providing quantitative, scale-independent measurements of the region segmented by one of the techniques described above. In this way, an absolute measurement of the region of interest in the image is possible, independent of the scale of the image, which could depend, for instance, on the distance of the objective to the camera.

\begin{figure}[!h]
\begin{subfigure}{0.23\textwidth}
\begin{center}
\includegraphics[height=3cm]{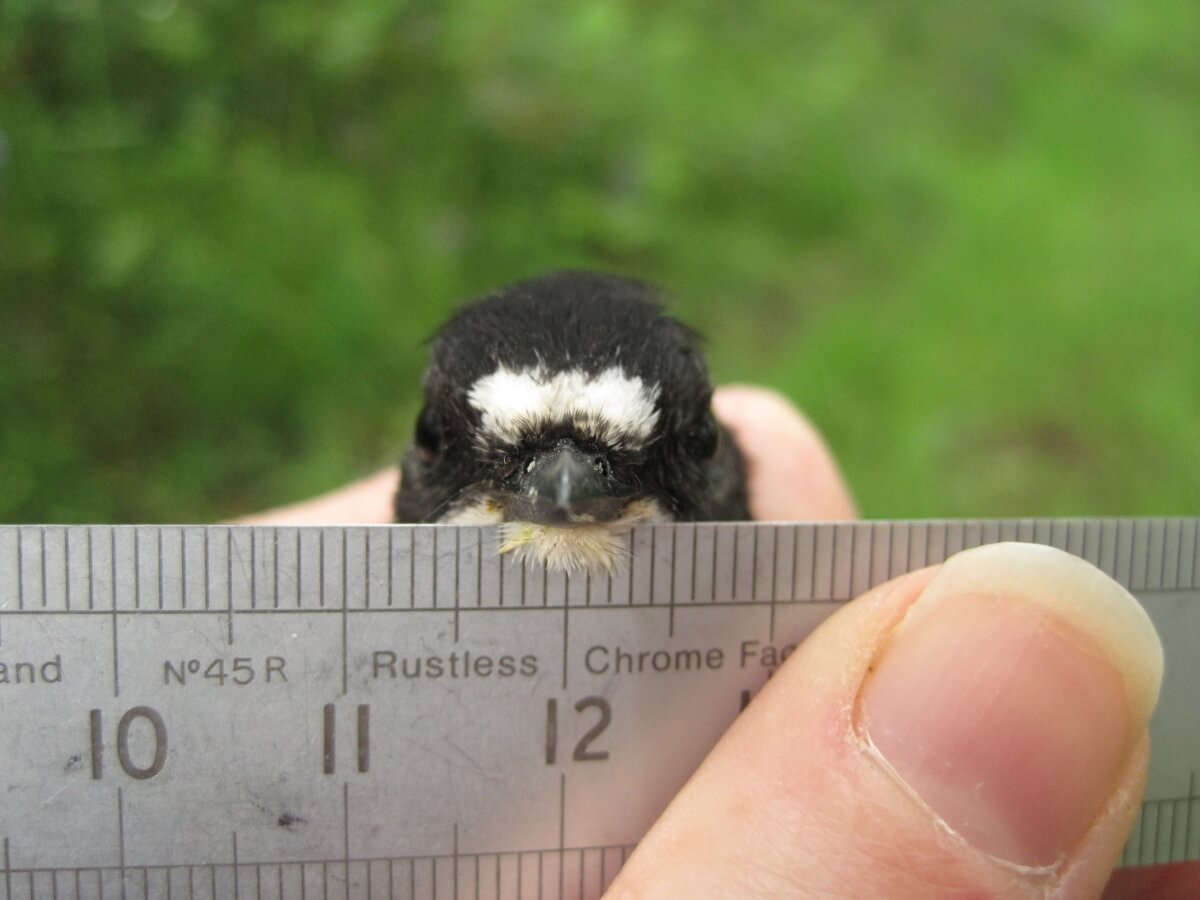}
\end{center}
\caption{}
\end{subfigure}
\begin{subfigure}{0.23\textwidth}
\begin{center}
\includegraphics[height=3cm]{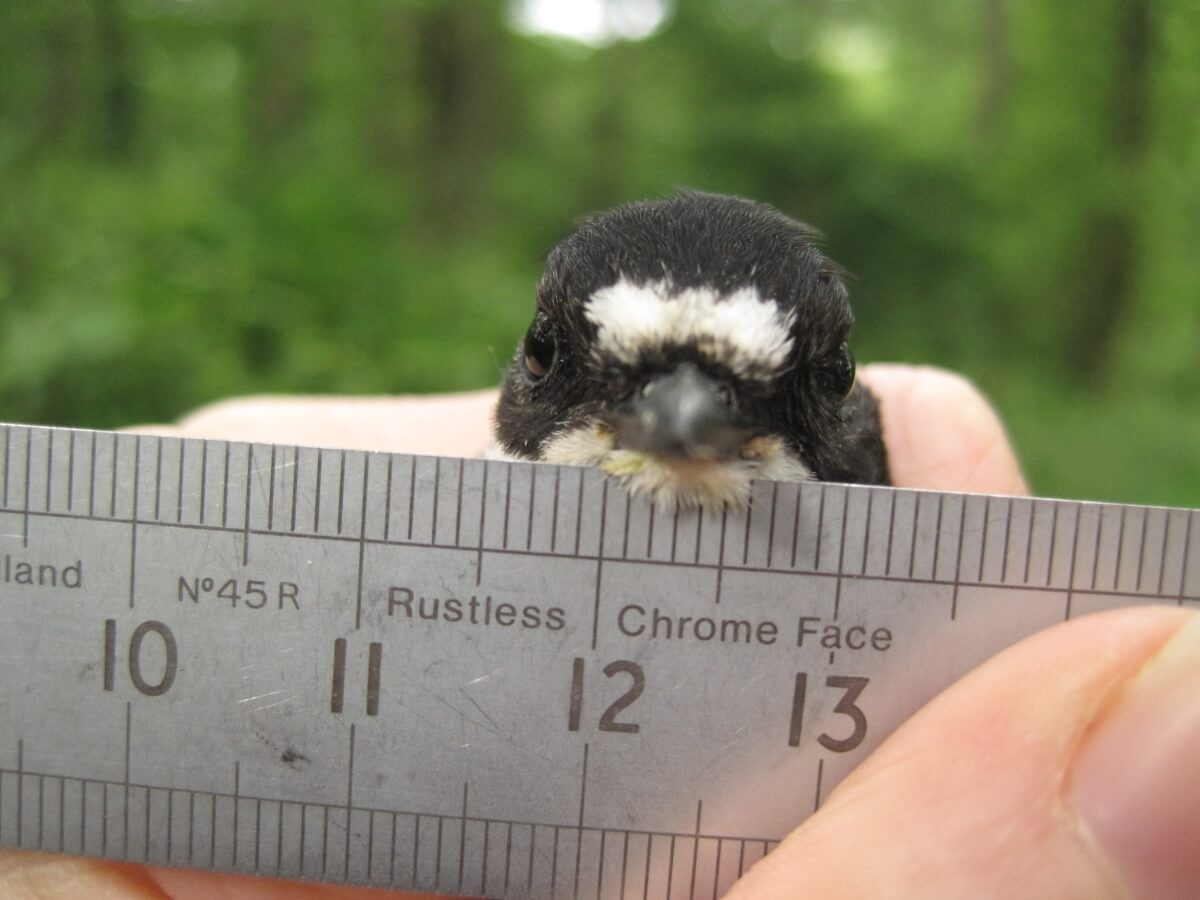}
\end{center}
\caption{}
\label{database:imb}
\end{subfigure}
\begin{subfigure}{0.23\textwidth}
\begin{center}
\includegraphics[height=3cm]{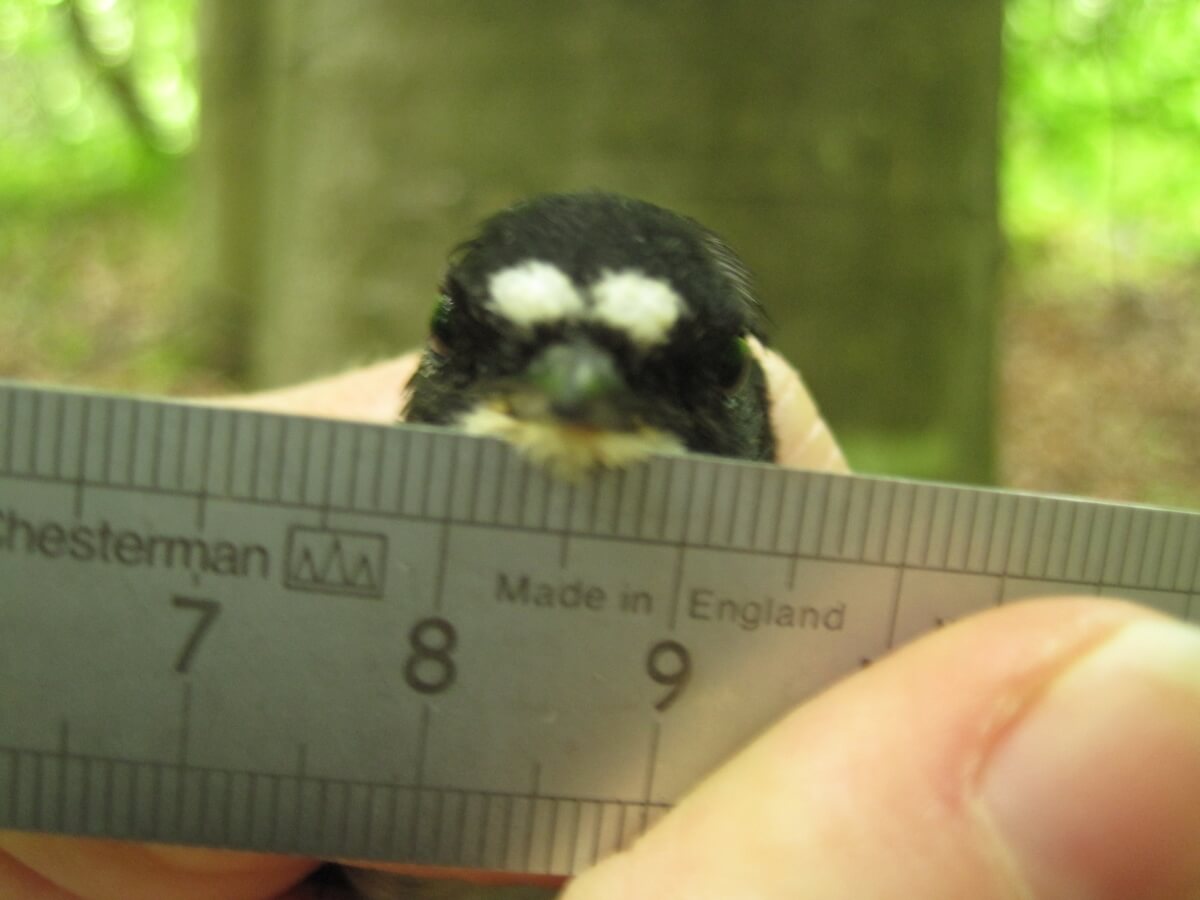}
\end{center}
\caption{}
\label{database:imc}
\end{subfigure}
\begin{subfigure}{0.23\textwidth}
\begin{center}
\includegraphics[height=3cm]{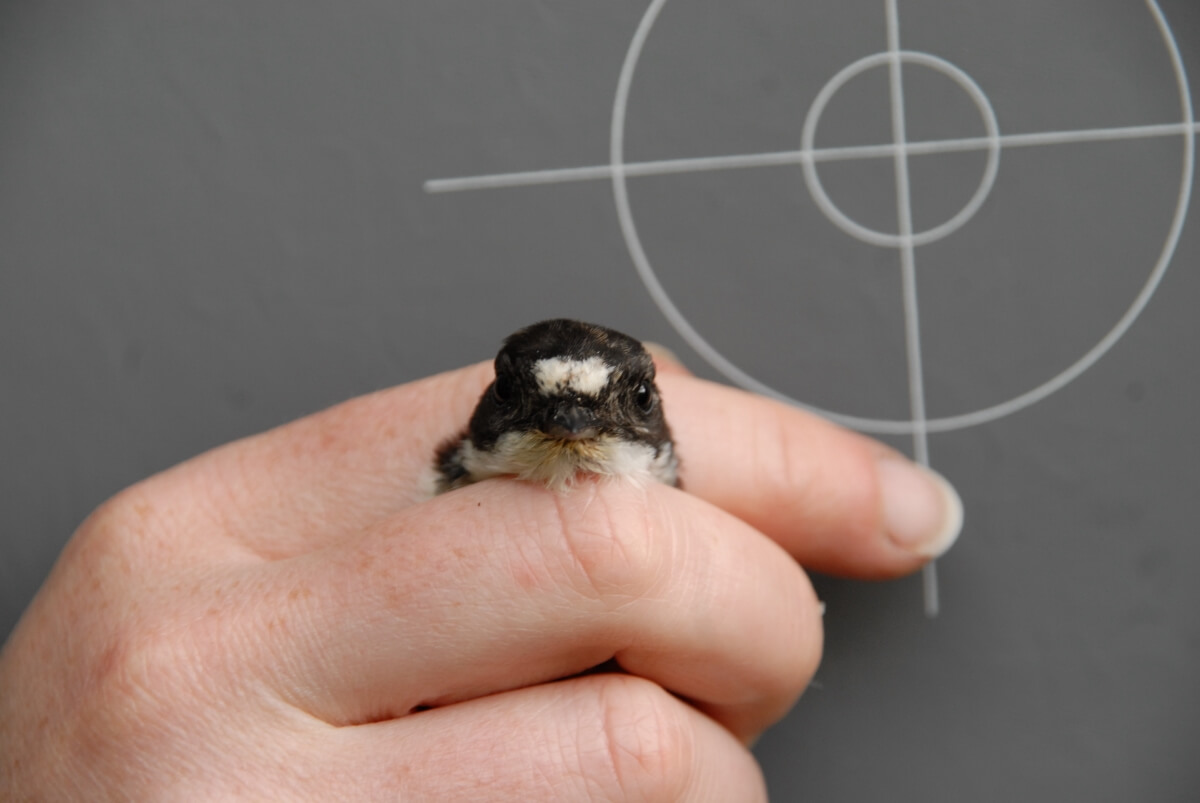}
\end{center}
\caption{}
\end{subfigure}
\caption{The blaze segmentation and measurement problem: pictures are taken at different distances, thus requiring a measurement tool.}
\label{fig:im_database}
\end{figure}

\begin{figure}[!h]
\begin{subfigure}{.5\textwidth}
\includegraphics[height=3cm]{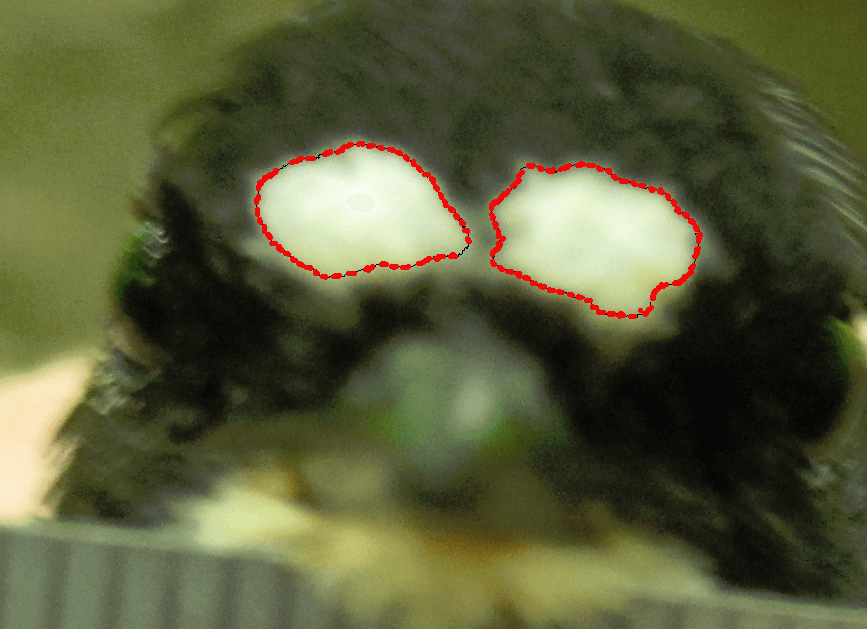}\ 
\includegraphics[height=3cm]{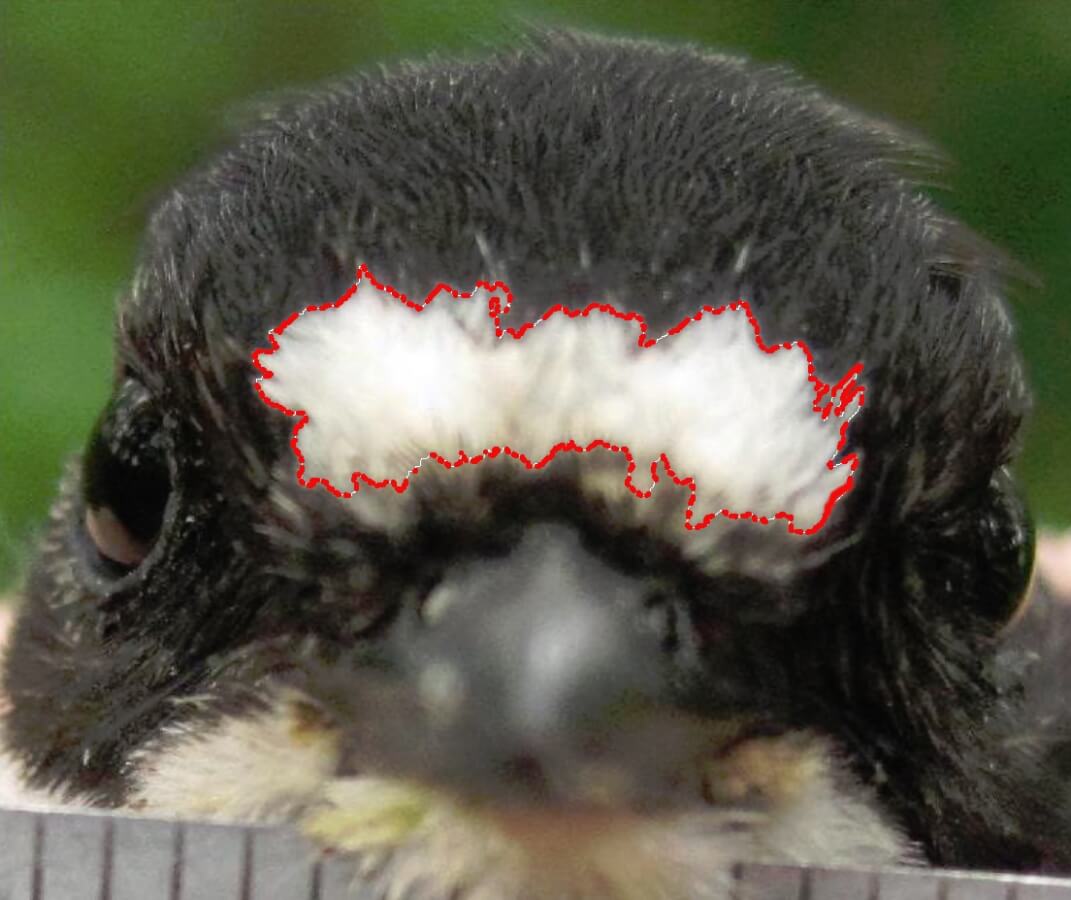}
\caption{Magic wand}
\label{fig:2a}
\end{subfigure}
\begin{subfigure}{.5\textwidth}
\includegraphics[height=3cm]{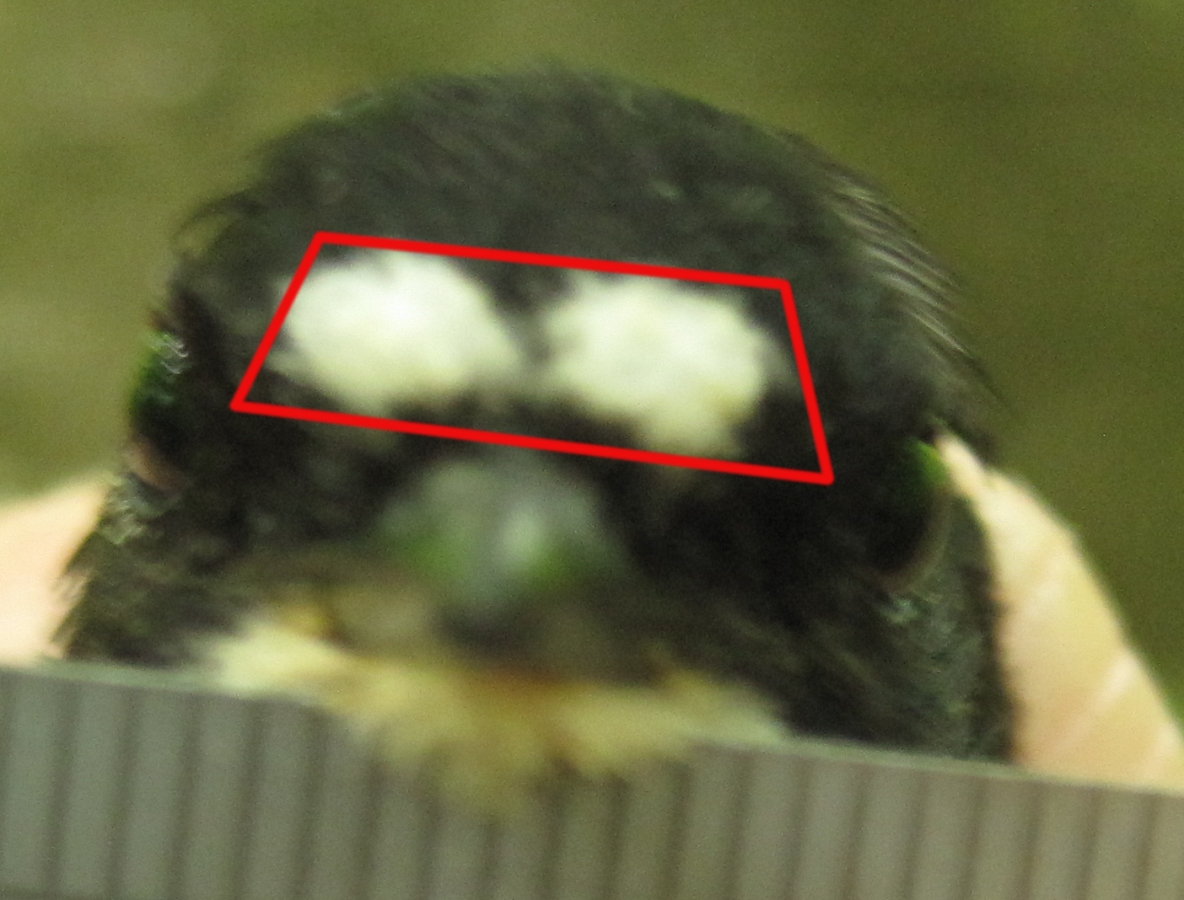}\
\includegraphics[height=3cm]{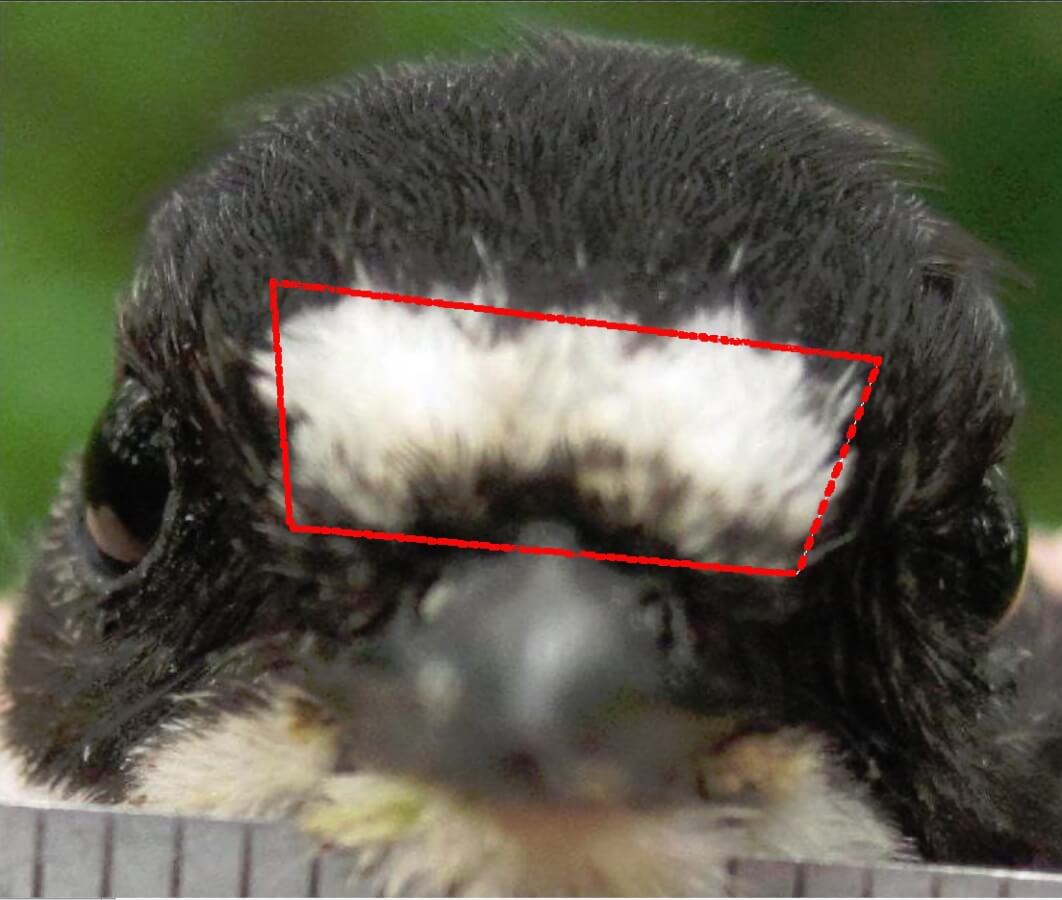}
\caption{Trapezium fitting}
\label{fig:2b}
\end{subfigure}
\caption{Flycatcher blaze segmentation of the images \ref{database:imb} and \ref{database:imc} obtained either by using the `magic wand' tool of the IMAGEJ software, similarly as described by Moreno \cite{morenovelando} or by trapezium fitting as suggested by Potti and Montalvo \cite{pottimontalvo}. In the first case the result is strongly user-dependent, in the second one the blaze area is overestimated.}
\label{fig:previous_methods_blaze}
\end{figure}

We demonstrate the use of our method in the context of real world applications in which segmentation and subsequent object measurement are crucial. Our main application is the measurement of the white forehead patch (blaze) of male pied flycatchers, which has been studied with regard to sexual selection in \cite{saetre}, see Figure \ref{fig:im_database}. The forehead patch is known to vary between individuals \cite{lundberg} and can be subject to both intra- \cite{jarvisto} and intersexual \cite{pottimontalvo} selection with pied flycatchers from Spain preferring males with large patches. Forehead patch size has been shown to signal male phenotypic quality through plasma oxidative stress and antioxidant capacity \cite{morenovelando}. However, in all studies to date the measurements of patches have been inconsistent and generally inaccurate. For example some studies have simply measured patch height \cite{daleslagsvold}, whereas Potti and Montalvo \cite{pottimontalvo} assumed the shape to be a trapezium with area equal to $0.5(B + b)H$, $B$ being the white patch width, $b$ the bill width and $H$ the height of the white patch. Morales et al. \cite{moralesmoreno} measured the length and breadth of the forehead patch with callipers to the nearest $0.01 mm$ and its size ($mm^2$) was calculated as the area of a rectangle. Other studies have measured the patches from photographs, e.g., J\"arvist\"o et al. \cite{jarvisto} Ruuskanen et al. \cite{ruuskanen} and Sirki\"a et al. \cite{sirkia} who photographed the forehead with a scale bar included in each picture, and measured the patch as the white area in $mm^2$ using IMAGEJ software \cite{Abramoff2004}. But none of these three papers provide methods of how the measurement was actually achieved, e.g., whether patches were delineated or roughly estimated with a simple shape. Most recently Moreno et al. \cite{morenovelando} analysed digital photos of forehead patches with Adobe PhotoShop CS version 11.0. relating the distance of $1 mm$ on the ruler to number of pixels, and used this to estimate length. Zooming to 400\% and using the paintbrush tool with 100\% hardness and 25\% spacing the authors delineate the patch and measure the area of the white areas on forehead. While this is the best measurement method to date, it still is subject to human measurement error and subjective assessment of patch boundaries. We report some segmentation results obtained by manual selection and polygon fitting in Figure \ref{fig:previous_methods_blaze}. In this manuscript we use a mathematically robust approach to segment the blaze independently to provide an accurate measurement of forehead patch area.

\begin{figure}[!h]
\begin{center}
\includegraphics[height=2.7cm]{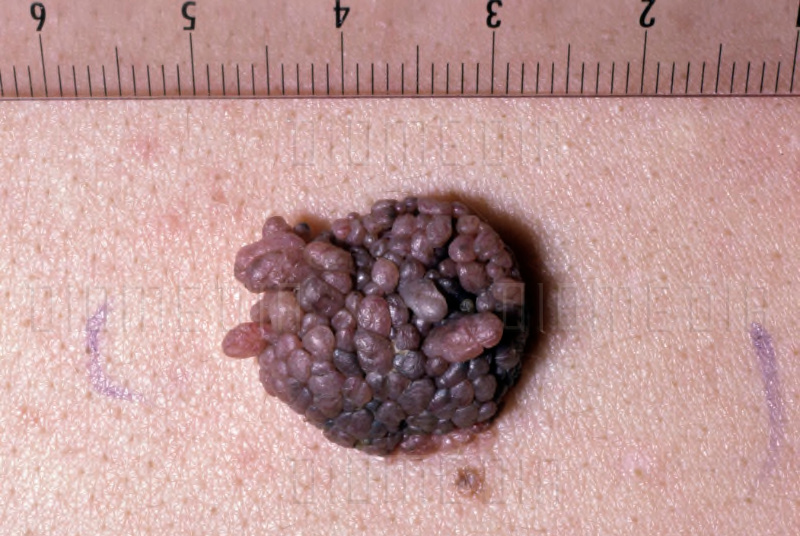}\hspace{0.7cm}
\includegraphics[height=2.7cm]{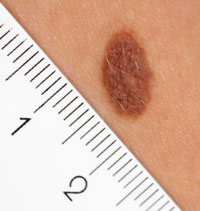}\hspace{0.7cm} 
\includegraphics[height=2.7cm]{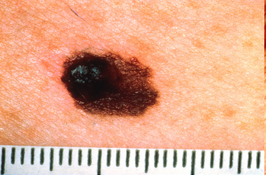}
\end{center}
\caption{The monitoring and measuring of moles is essential for the early diagnosis of melanoma. Normally, due to their small size, they can be measured by juxtaposing a small ruler with them.}
\label{fig:skinmoles}
\end{figure}

A similar challenge can be encountered in medical applications monitoring and quantifying the evolution of skin moles for early diagnosis of melanoma (skin cancer). A normally user-dependent measurement of the mole is performed using a ruler located next to it. A picture is then taken and used for future comparisons and follow-up, see Figure \ref{fig:skinmoles} and compare \cite{cavalcanti,abbas} for previous attempts of automatic detection of melanomas. For such an application, a systematic quantitative analysis is also required 
\footnote{Mole images from
\url{http://www.medicalprotection.org/uk/practice-matters-issue-3/skin-lesion-photography}, 
\textcopyright{Chassenet/Science Photo Library} ,  \\ \url{http://en.wikipedia.org/wiki/Melanoma} (public domain), \\
\url{http://www.diomedia.com/stock-photo-close-up-of-a-papillomatous-dermal-nevus-mole-a-raised-pigmented-skin-lesion-that-results-from-a-proliferation-of-benign-melanocytes-c-cid-image14515019.html}, \textcopyright{Phototake RM/ ISM}}. 

In several other applications the task of measuring objects directly from the image is encountered. These include zoological and behavioural studies arising in the  animal world where detecting size, shape and possible symmetries of specific distinctive animal features can be useful, as well as, for instance, in archaeological digs where the measurement of finds is important for comparisons and classification \cite{herrmann}. 

\paragraph{Outline of the method.} We consider the image as a graph whose vertices are the image pixels. Similarity between pixels in terms of colour or texture features is modelled by a weight function defined on the set of vertices. Our method runs as follows. Firstly, using examples provided by the user (dictionaries) as well as matrix completion and operator splitting techniques, the segmentation of the region of interest is performed. In the graph framework, this corresponds to cluster together pixels having similar features. This is obtained by minimising on the graph the Ginzburg-Landau functional typically used in the continuum setting to describe diffuse interface problems. In order to provide quantitative measurements of the segmented region, a second detection step is then performed. The detection here aims to identify the distinctive geometrical features of the measurement tool (such as line alignment for rulers or circularity for circles) to get the scale on the measurement tool considered. The segmented region of interest can now be measured by simple comparisons and quantitative measurements such as perimeter and area can be provided.

\paragraph{Contribution} We propose a self-contained programme combining automated detection and subsequent size measurement of objects in images where a measurement tool is present. Our approach is based on two powerful image analysis techniques in the literature: a graph segmentation approach which uses a discretised Ginzburg-Landau energy \cite{BerFlen} for the detection of the object of interest and the Hough transform \cite{hough} for detecting the scale of the measurement tool. While these methods are state of the art, their combination for measuring object size in images proposed in this paper is new. Moreover, to our knowledge there is only little contribution in the literature that broach the issue of how the graph segmentation approach as well as the Hough transform are applied to specific problems \cite{BerFlen1,BerFlen2,grah}. Indeed, here we present these methodologies in detail, especially discussing important aspects in their practical implementation, and demonstrate the robust applicability of our programme for measuring the size of objects, showcasing its performance on several examples arising in zoology, medicine and archaeology. Namely, we first apply our combined model for the measurement of the blaze on the forehead of male pied flycatchers, for which we run a statistical analysis on the accuracy and predicted error in the measurement on a database of thirty images. State-of-the-art methods for such a task typically require the user to fit polygons inside or outside the blaze  \cite{pottimontalvo} or to segment the blaze by hand \cite{morenovelando}. Similarly, the scale on the measurement tool is typically read from the image by manually measuring it on the ruler. With respect to medical applications, we apply our combined method for the segmentation and measurement of melanomas. Although efficient segmentation methods for automatised melanoma detection already exist in literature (see, e.g., \cite{cavalcanti,abbas}), up to the knowledge of the authors no previous methods providing their measurement by detecting the scale on the the ruler placed next to them (see Figure \ref{fig:skinmoles}) exist. Conversely, in the case of archaeological applications, some models for the automatic detection of the measurement tool in the image exist \cite{herrmann} but no automatic methods are proposed for the segmentation of the region of interests. \red{A free release of the MATLAB code used to compute the results will be made available after the zoological analysis of the pied flycatcher's data based on our segmentation and measurement has been completed, \cite{ourpreprint}.}

\paragraph{Organisation of the paper.}  In Section \ref{sec:graph} we present the mathematical ingredients used for the design of the graph based segmentation technique used in \cite{BerFlen,BerFlen1,BerFlen2,MerkujevBert2015}. They come from two different worlds:  the framework of diffusion PDEs used for modelling phase transition/separation problems (see Section \ref{subsec:GL}) and graph theory and clustering, see Section \ref{subsec:graph}. In view of a detailed numerical explanation, we also recall a splitting technique and a popular matrix completion technique used in our problem to overcome the computational costs. \red{In Section \ref{sec:hough_use} we explain how the geometrical Hough transform is used to detect the scale in an image}. Finally, Section \ref{sec:res} contains the numerical results obtained with our combined method applied to the problems described above. \red{For completeness, we give some details on the Nystr\"{o}m matrix completion technique in Appendix \ref{nystrom:app} and a review of the Hough transform for line and circle detection in Appendix \ref{sec:hough}.}

\section{Image segmentation as graph clustering} \label{sec:graph}

We present in this section the mathematical background for the design of the Ginzburg-Landau (GL) graph segmentation algorithm introduced in \cite{BerFlen}. There, the image segmentation problem is rephrased as a minimisation problem on a graph defined by features computed from the image. Compared to the methods above, the graph framework allows for more freedom in terms of the possible features used to describe the image, such as texture.

\subsection{The Ginzburg-Landau functional as approximation of TV} \label{subsec:GL}

In the following, we recall the main properties of the original continuum version of the GL functional explaining its importance in the context of image segmentation problems as well as the main concepts of graph theory which will be used for the segmentation modelling.

Several physical problems modelling phase transition and phase separation phenomena are built around the well-known GL functional:
\begin{equation} \label{def:GL_func}
GL(u) := \frac{\varepsilon}{2} \int_\Omega |\nabla u(x)| ^2\ dx +\frac{1}{\varepsilon}\int_\Omega W(u(x))\ dx .
\end{equation}
The functional above is defined in the continuous setting. Here, $\Omega$ represents a open subset of $\R^d, d=2,3$, $u:\Omega\to\R$ is the density of a two-phase material and $W(u)$ is a double-well potential, e.g. $W(u)=\frac{1}{4}(u^2-1)^2$.  The two wells $\pm 1$ of $W$ correspond to the two phases of the material. The parameter $\varepsilon>0$ is the spatial scale. Variational models built around this functional are also referred to as \emph{diffuse interface} models because of the interface appearing between the two regions containing the phases (i.e. the two wells of $W$) due to the competition between the two terms of the functional \eqref{def:GL_func}. Nonetheless, some smoothness preventing $u$ from having jumps between the two regions is ensured by the first regularisation term.

The use of the GL functional has become very popular in image processing due to its connections with the total variation (TV) seminorm. In \cite{modicamortola1,modicamortola2}, for instance, $\Gamma$-convergence properties of \eqref{def:GL_func} to the TV functional are shown. Thus, the GL functional is very often used as a quadratic approximation of total variation. Fast numerical schemes relying on these connections have been designed for many imaging problems, thus overcoming the issues related to nonsmooth TV minimisation \cite{ambrosio,esedoglu,chanvese}. In image processing, the functional considered often is of the form
\begin{equation}  \label{func:min}
E(u):=GL(u)+ \lambda\ \phi(u,u_0),
\end{equation}
where $\phi(u,u_0)$ is a fidelity term measuring the distance of the reconstructed image $u$ to the given image $u_0$. Depending on the application, different data fidelities are employed. Typically, they are related to statistical and physical assumptions of the model considered. Standard examples of fidelity terms are $\phi(u,u_0)=\norm{u-u_0}^d_{L^d(\Omega)},~ d=1,2$. The parameter $\lambda>0$ determines the influence of the data fit compared to the regularisation. Taking the $L^2$ gradient descent of \eqref{func:min} we get the following evolutionary PDE, known in the literature as the Allen-Cahn equation \cite{alcahn} with an additional forcing term due to the fidelity $\phi$:
\begin{equation} \label{grad_desc}
u_t=-\frac{\delta GL}{\delta u} - \lambda \frac{\delta \phi}{\delta u} = \varepsilon \Delta u -\frac{1}{\varepsilon} W'(u) -\lambda \frac{\delta \phi}{\delta u}.
\end{equation}
Steady states of equation \eqref{grad_desc} are the stationary points of the energy $E$ in \eqref{func:min}. Note that $E$ is not convex so uniqueness is not guaranteed and, consequently, the long time behaviour for solutions of \eqref{grad_desc} will depend  on the initial condition. The linear diffusion term weighted by $\varepsilon$ appearing in \eqref{grad_desc} allows for fast solvers using for instance the Fast Fourier Transform (FFT) which translates the Laplace operator into a multiplication operator on the Fourier modes.


\subsection{Towards the modelling: the graph framework}  \label{subsec:graph}

In the following, we rely on the method presented in \cite{BerFlen,MerkujevBert2015} for high-dimensional data classification on graphs which has been applied to several imaging problems \cite{BerFlen1,BerFlen2}, showing good performance and robustness. We consider the problem of binary image segmentation where we want to partition a given image into two components where each component is a set of pixels (also called a cluster, or a class) and represents a certain object or group of objects. Typically, some \emph{a priori} information describing the object(s) we want to extract is given and serves as initial input for the segmentation algorithm.  For image labelling, in \cite{BerFlen} two images are taken as input: the first one has been manually segmented in two classes and the objective is to automatically segment the second image using the information provided by the segmentation of the first one.

\smallskip

We revise in the following the main ingredients of the model considered and start from a quick review of concepts in graph theory. \red{We represent a rectangular image with  $S:=N\times M$ pixels by the set $I:= \{x=(x_1,x_2)\in \Z^2: 0\leq x_1 \leq N-1 \text{ and } 0\leq x_2 \leq M-1\}$}. For each $x\in I$, we define the image neighbourhood of $x$ as the set 
$$\mathcal{N}(x):=\left\{y\in I: |x_1-y_1|\leq \tau \text{ and }|x_2-y_2|\leq \tau\right\},$$
with $\tau\in\mathbb{N}$ fixed, i.e. $\mathcal{N}(x)$ contains the pixels in a $(2\tau+1) \times (2\tau+1)$ sized square centred at $x$. For some appropriate $K\in \N$, we associate to every pixel $x\in I$ a vector $z\in\R^K$ encoding selected characteristics of the neighbourhood $\mathcal{N}(x)$. These characteristics are related to the grey or RGB (red, green, blue) intensity values as well as the texture features of the neighbourhood. In Section~\ref{subsec:algGL}, we will explain in more detail our feature vector construction. The map $\psi: I\to\R^K, x \mapsto z$ is called the feature function. For constructing the feature vectors \red{in Section \ref{subsec:algGL},} it will be useful to associate a neighbourhood vector $\nu(x) := (x_j)_{j\in \mathcal{N}(x)} \in I^{(2\tau+1)\times (2\tau+1)}$ to each neighbourhood, such that the ordering of the $x_j$ in $\nu(x)$ is consistent between pixels $x$, e.g., order the pixels from each square $\mathcal{N}(x)$ from left to right and top to bottom. The specific choice of ordering is not important, as long as it is consistent for each pixel neighbourhood. 

Next we construct a simple weighted undirected graph $G=(V,E,w)$ whose vertices correspond to the pixels in $I$ and with edges whose weights depend on the feature function $\psi$. Let $V$ be a vertex set of cardinality $S$. To emphasize that each vertex in $V$ corresponds to exactly one pixel in $I$, we will label the vertex corresponding to $x\in I$ by $x$ as well. Let $w:V\times V\to\R$ be a symmetric and nonnegative function, i.e. for each $x_i,x_j\in V$
\begin{equation} \label{defin:weightfunct}
w(x_i,x_j)=w(x_j,x_i),\qquad  w(x_i,x_j)\geq 0.
\end{equation}
We define the edge set $E$ as the collection of all undirected edges connecting nodes $x_i$ and $x_j$ for which $w(x_i,w_j) >0$ \cite{chung}. The function $w$ restricted to $E\subset V\times V$ is then a positive edge weight function.

In our applications we define $w$ as 
\begin{equation*}
w(x_i,x_j):=\hat{w}(\psi(x_i),\psi(x_j))=\hat{w}(z_i,z_j),
\end{equation*}
where $\hat{w}:\R^K\times\R^K\to\R$ is a given function and $\psi$ is the feature function.

In operator form, the \textbf{weight matrix} $W\in\R^{S\times S}$ is the nonnegative symmetric matrix whose elements are $w_{i,j}=w(x_i,x_j)$.
In the following, we will not distinguish between the two functions $w$ and $\hat{w}$ and, with a little abuse of notation, we will write $w(z_i,z_j)$ for $\hat{w}(z_i,z_j)$.

\begin{remark}
Weight functions express the similarities between vertices and will be used in the following  to partition $V$ into clusters such that the sum of the edge weights between the clusters is small.  There are many different mathematical approaches to attempt this partitioning. When formulated as a balanced cut minimisation, the problem is NP-complete \cite{luxburg}, which inspired relaxations which are more amenable to computational approaches, 
many of which are closely related to spectral graph theory \cite{shimalik}. We refer the reader to \cite{chung} for a monograph on the topic. The method we use in this paper can be understood (at least in spirit, if not technically, \cite{vGB,vanGennip}) as a nonlinear extension of the linear relaxed problems.
\end{remark}

To solve the segmentation problem, we minimise a {\it discrete GL functional} (which is formulated in the graph setting, instead of the continuum setting), via a gradient descent method similar to the one described in Section \ref{subsec:GL}. In particular, in this setting the Laplacian in \eqref{grad_desc} will be a (negative) normalised graph Laplacian. We will use the spectral decomposition of $u$ with respect to the eigenfunctions of this Laplacian. In Section \ref{subsec:nystrom} we discuss the Nystr\"om method, which allows us to quickly compute this decomposition, but first we introduce the graph Laplacian and graph GL functional.

%
%

\paragraph{The discrete operators.}
We start from the definition of the differential operators in the graph framework.

For each vertex $x\in V$, we define the \textbf{degree} of $x$, 
\begin{equation*}
d:V\to\R,\quad d(x):=\sum_{y\in V} w(x,y).
\end{equation*}
In operator form, the diagonal \textbf{degree matrix} $D\in\R^{S\times S}$ is defined to have diagonal elements $d_{i,i}=d(x_i)$.

A subset $A$ of the vertex set $V$ is \textbf{connected} if any two vertices in $A$ can be connected by a path (i.e. a sequence of vertices such that subsequent vertices are connected by an edge in $E$) such that all the vertices of the path are in $A$.
A finite family of sets $A_1,\ldots,A_t$ is called a \textbf{partition} of the graph if $A_i\cap A_j=\emptyset$ for  $i\neq j$ and $\bigcup_i A_i=V$.

\medskip

We now have all the ingredients to define the \textbf{graph Laplacian}. Denoting by $\mathcal{V}$ the space of all the functions $V\to\R$, the graph Laplacian is the operator $L:\mathcal{V}\to\mathcal{V}$ such that:
\begin{equation}  \label{def:graph_laplofafunc}
Lu(x)=\sum_{y\in V} w(x,y)(u(x)-u(y)),\quad x\in V.
\end{equation}
We are considering a finite graph of size $S$, so real valued functions can be identified as vectors in $\R^S$. We can then write the graph Laplacian in matrix form as $L=D-W$ or element-wise as:
\begin{align} \label{def:graph_lapl}
&L(x,y):=\left\{\begin{aligned}
d(x),\quad&\text{ if }x=y,\\
-w(x,y),\quad&\text{ otherwise.}
\end{aligned}\right.
\end{align}
It is worth mentioning (see Remark~\ref{rem:spectrum} below) that this graph Laplacian is a positive semidefinite operator. Note that by convention the sign of the discrete Laplacian is opposite to that of the (negative semidefinite) continuum Laplacian.
The associated quadratic form of $L$ is
\begin{equation}   \label{def:quadr}
Q(u,Lu):=\frac{1}{2}\sum_{x,y\in V} w(x,y)\left(u(x)-u(y)\right)^2.
\end{equation}
The quadratic form $Q$ can be interpreted as the energy whose optimality condition corresponds to the vanishing of the graph Laplacian in \eqref{def:graph_lapl}.

\begin{remark}\label{rem:spectrum}
The operator $L$ has $S$ non-negative real-valued eigenvalues $\left\{\lambda_i\right\}_{i=1}^S$ which satisfy: $0=\lambda_1\leq \lambda_2\leq\cdots\leq \lambda_S$. The eigenvector corresponding to $\lambda_1$ is the constant $S$-dimensional vector \textbf{1}$_S$, see \cite{luxburg}.
\end{remark}

The operator in \eqref{def:graph_laplofafunc}-\eqref{def:graph_lapl} is not the only graph Laplacian appearing in the literature. To set it apart from others, it is also referred to as the unnormalised or combinatorial graph Laplacian. Such operator can be related to the standard continuous differential one through nonlocal calculus \cite{Gilboa2}. More precisely, the eigenvectors of $L$ converge to the eigenvectors of the standard Laplacian, but in the large sample size limit a proper scaling of $L$ is needed in order to guarantee stability of convergence to the continuum operator \cite{BerFlen,BerFlen2}. Hence, we consider in the following the normalisation of $L$ given by the symmetric graph Laplacian
\begin{equation}  \label{symmgraphlapl}
L_s:=D^{-1/2}LD^{-1/2}=I-D^{-1/2}WD^{-1/2}.
\end{equation}
Clearly, the matrix $L_s$ is symmetric. Other normalisations of $L$ are possible, such as the random walk graph Laplacian (see \cite{chung,luxburg,vanGennip}).

\medskip

In \cite[Section 5]{shimalik} a quick review on the connections between the use of the symmetric graph Laplacian \eqref{symmgraphlapl} and spectral graph theory is given. Computing the eigenvalues of the normalised symmetric Laplacian corresponds to the computation of the generalised eigenvalues used to compute normalised graph cuts in a way that the standard graph Laplacian may fail to do, compare \cite{chung}. Typically, spectral clustering algorithms for binary segmentation base the partition of a connected graph on the eigenvector corresponding to the second eigenvalue of the normalised Laplacian, using, for example, $k$-means. For further details and a comparison with other methods we refer the reader to \cite{shimalik} and to \cite[Section 2.3]{BerFlen} where a detailed explanation on the importance of the normalisation of the Laplacian is given.

\paragraph{The discrete GL functional.}
Recalling \eqref{def:GL_func}-\eqref{func:min} and \eqref{def:quadr}, we define the discrete GL functional\footnote{`Discrete GL functional with a data fidelity term' would be a more accurate name, but we opt for brevity here.} as
\begin{equation}   \label{def:GL_discr}
GL_d(u): =\frac{\varepsilon}{2}\ Q(u,L_s u) + \frac{1}{\varepsilon} \sum_{x\in V} W(u(x))+ \sum_{x\in V} \frac{\chi(x)}{2}(u(x)-u_0(x))^2. \notag
\end{equation}
Here $u_0$ represents known training data provided by the user. As before, $W(u(x))=\frac{1}{4}(u^2(x)-1)^2$ is the double-well potential. The function $\chi: V\to \{0,1\}$ is the characteristic function of the subset of labelled vertices $V_{lab}\subset V$, i.e. $\chi=1$ on $V_{lab}$ and $\chi=0$ on $V_{unlab}:=V_{lab}^c$. Hence, the corresponding fidelity term enforces the fitting between $u$ and $u_0$ in correspondence to the the known labels on the set $V_{lab}$, while the labelling for the pixels in $V_{unlab}$ is driven by the first two regularising terms in \eqref{def:GL_discr}.

The corresponding $\ell^2$ gradient flow for  \eqref{def:GL_discr} reads
\begin{equation*}
u_t=-\varepsilon~L_s u -\frac{1}{\varepsilon}  \sum_{x\in V} (u^3(x)-u(x))- \sum_{x\in V} \chi(x)(u(x)-u_0(x)).
\end{equation*}

The idea is to design a semi-supervised learning (SSL) approach where \emph{a priori} information for the set $V_{lab}$ (i.e. cluster labels) is used to label the points in the set $V_{unlab}$. The comparison uses the weight function defined in \eqref{defin:weightfunct} to build the graph by comparing the feature vectors at each point.

\begin{remark}[The weight function]
As pointed out in \cite[Section 2.5]{BerFlen}, the main criteria driving the choice of the weight function are the desired outcome and the computational efforts required to diagonalise the corresponding matrix $W$. A common weight function is the Gaussian function, which, for $x,y\in V$ reads
\begin{equation}   \label{similarity_func}
w(x,y)=\exp (-\|\psi(x)-\psi(y)\|^2/\sigma^2),\quad\sigma>0.
\end{equation}
Note that this function is symmetric: $w(x,y)=w(y,x)$.
\end{remark}

Several approaches to SSL using graph theory have been considered in literature, compare \cite{Coifman,Gilboa2}.
The approach presented here adapts fast algorithms available for the efficient minimisation of the continuous GL functional to the minimisation of the discrete one in \eqref{def:GL_discr} . In particular, to overcome the high computational costs, we present in the following an operator splitting scheme and a matrix completion technique applied to our problem.

\subsection{Convex splitting}  \label{subsec:convexsplit}

Splitting methods are used in the study of PDEs. 
Here, we focus on convex splitting, 
which is used to numerically solve problems with a general gradient flow structure. \red{Decomposing $GL_d$ as
$$
GL_d=GL_{1,d}-GL_{2,d}
$$  where both $GL_{1,d}$ and $GL_{2,d}$ are convex
and denoting by $U_n$ the spatially discretisation of $u(\cdot,n\Delta t)$, $\Delta t>0$, $n\geq0$, a semi-implicit discretisation for the steepest descent of $GL_d$ reads
\begin{equation}   \label{convexsplitphilosophy}
U_{n+1}-U_n=-\Delta t (\nabla_V GL_{d,1}(U_{n+1})-\nabla_V GL_{d,2}(U_n)),
\end{equation}
where $\nabla_V$ indicates formally the Fr\'echet derivative with respect to the metric in a Banach space $V$.}
The advantage of the convex splitting consists in treating the convex
part implicitly in time and the concave part explicitly. Typically, nonlinearities are considered in the explicit part of the splitting and their instability is balanced by the effect of the implicit terms.


The terms $GL_{d,1}$ and $GL_{d,2}$ in \eqref{convexsplitphilosophy} read in our case (cf. \cite[Section 3.1]{BerFlen})
\begin{subequations} \label{GLdconvexsplit}
\begin{equation} \label{convexsplit1}
GL_{d,1}(u):=\frac{\varepsilon}{2}~ Q(u,L_s u) +\frac{C}{2}\sum_{x\in V}  u^2(x),
\end{equation}
\begin{equation}\label{convexsplit2}
GL_{d,2}(u):= -\frac{1}{4\varepsilon}\sum_{x\in V} (u^2(x)-1)^2 +\frac{C}{2}\sum_{x\in V}  u^2(x) - \sum_{x\in V}\frac{\chi(x)}{2} (u(x)-u_0(x))^2, 
\end{equation}
\end{subequations}
where the constant $C>0$ has to be chosen large enough such that $GL_{d,2}$ is convex for $u$ around the wells of $W$. The differential operator contained in the implicit component of the splitting, $GL_{d,1}$, is the symmetric graph Laplacian, which can be diagonalised quickly and inverted using Fourier transform methods. In \cite[Section 3.1]{BerFlen}, more details of the splitting are presented.
\red{Writing out in detail the time-discretised scheme \eqref{convexsplitphilosophy}, we get, for every $n\geq 1$
\begin{multline}   \label{GLconvexsplit}
U_{n+1}(x)-U_n(x)= -\Delta t\left( \varepsilon~L_s(U_{n+1}(x)) + C U_{n+1}(x)\right) -\Delta t\left(- \frac{1}{\varepsilon} \left(U^3_{n}(x) - U_{n}(x) \right) +  C ~U_{n}(x)\right.\\  - \chi(x)\left(U_{n}(x) -U_0\right) \Big),\quad x\in V.
\end{multline}
Here, $U_0$ denotes the training data, i.e. the known labels $-1$ and $1$ assigned by the user to nodes in the subset $V_{lab}\subset V$.  In our numerical experiments
we initialised the time-stepping \eqref{GLconvexsplit} by taking 
\begin{equation}   \label{initialisation:GL}
U_1(x)=\left\{
\begin{aligned}
& U_0(x),\quad &\text{if }x\in V_{lab},\\
& 0,\quad &\text{if }x\in V^C_{lab}.
\end{aligned} \right.
\end{equation}
}

\paragraph{Towards the numerical realisation.} The numerical strategy we intend to use is based on the following steps (see Section \ref{subsec:algGL} for more details):

 \begin{itemize}
 \item At each time step $n\Delta t, n\geq 1$, consider at every point the spectral decomposition of $U_n$ with respect to the eigenvectors $v_k$  of the operator $L_s$ as
 \begin{equation}  \label{spectrdecun}
 U_n(x)=\sum_{k} \alpha_n^k(x) v_k(x),\quad x\in V
 \end{equation}
 with coefficients $\alpha_n$.
Similarly, use spectral decomposition in the $\left\{v_k\right\}$ basis for the other nonlinear quantities appearing in \eqref{GLconvexsplit}.
 \item Having fixed the basis of eigenfunctions, the numerical approximation in the next time step $U_{n+1}$ is computed by determining the new coefficients $\alpha^k_{n+1}$ in \eqref{spectrdecun} for every $k$ through convex splitting \eqref{GLconvexsplit}. 
  \end{itemize}

The only possible bottleneck of this strategy is the computation of the eigenvectors $v_k$ of the operator $L_s$, which, in practice, can be computationally costly for large and non-sparse matrices $W$. 
To mitigate this potential problem, we use the Nystr\"om extension (Section~\ref{subsec:nystrom}).

\subsection{Matrix completion via Nystr{\"o}m extension}  \label{subsec:nystrom}

Following the detailed discussion in \cite[Section 3.2]{BerFlen}, we present here the Nystr{\"o}m technique for matrix completion \cite{nystrom} used in previous works by the graph theory community \cite{fowlkes,belongie} and applied later to several imaging problems \cite{naeini,MerkujevKostic2013,MerkujevSunu2014}. In our problem, the Nystr{\"o}m extension is used to find an approximation of the eigenvectors $v_k$ of the operator $L_s$. We will freely switch between the representation of eigenvectors (or eigenfunctions) as a real-valued functions on the vertex set $V$ and as a vectors in $\R^S$.

\smallskip

Consider a fully connected graph with vertices $V$  and the set of corresponding feature vectors $\psi(V)=\{z_i\}_{i=1}^S$. A vector $v$ is an eigenvector of the operator $L_s$ in \eqref{symmgraphlapl} with eigenvalue $\lambda$ if and only if $v$ is an eigenvector of the operator $D^{-1/2}WD^{-1/2}$ with eigenvalue $1-\lambda$, since
\begin{align}\label{eigenvectorsLs}
& L_s v = v - D^{-1/2}WD^{-1/2} v = \lambda v \quad \Longleftrightarrow\\ 
& D^{-1/2}WD^{-1/2} v = (1-\lambda) v. \notag
\end{align}
Thus, finding the spectral decomposition of $L_s$ boils down to diagonalising the operator $D^{-1/2}WD^{-1/2}$. This is not obviously easier, as the matrix $W$, despite being nonnegative and symmetric, may be large and non-sparse, so the computation of its spectrum may be computationally hard. Here, however, we take advantage of the Nystr{\"o}m  extension. Given the eigenvalue problem
\begin{align} \label{eig:eq}
& \text{find }\theta\in\R\text{ and }v: V\to \R, v\neq0\quad\text{ s. t. } \\
& \sum_{x\in V} w(x,y)~v(x) = \theta v(y), \notag
\end{align}
for every point $y\in V$, we approximate the sum on the left hand side using a standard quadrature rule where the interpolation points are chosen by randomly selecting a subset of $L$ points from the set $V$ and the interpolation weights are chosen correspondingly.  The Nystr{\"o}m extension for \eqref{eig:eq} then approximates  \eqref{eig:eq} by
\begin{equation} \label{nystrom1}
\sum_{i=1}^{L} w(y,x_i) v(x_i)\approx\sum_{x\in V} w(y,x) v(x)= \theta v(y),
\end{equation}
where $X:=\{x_i\}_{i=1}^L\subset V$ is a set of randomly chosen vertices.
The set $X$ defines a partition of $V$ into $X$ and $Y:=X^c$. In \eqref{nystrom1} we approximate the value $v(y)$, for an eigenvector $v$ of $W$ and $y\in Y$, only knowing the values $v(x_i),~ i=1,\ldots,L$, by solving the linear problem
\begin{equation} \label{nystrom2}
\sum_{i=1}^{L} w(y,x_i) v(x_i)=\theta v(y).
 \end{equation}
With this method we can approximate the values of an eigenvector $v$ of $W$, corresponding to the eigenvalue $\theta$, in the \emph{whole} set of points $V$ using its values in the subset $X$ and solving the interpolated eigenvalue equation above. Generally, this is not as immediate as it sounds since the eigenvectors of $W$ are not known in advance, however, by choosing $y=x_j$, $j=1, \ldots, L$, in \eqref{nystrom2}, we find an eigenvalue problem for the known matrix with entries $w(x_j,x_i)$, which is a much smaller matrix than the full matrix $W$:
\begin{equation} \label{nystrom3}
\sum_{i=1}^{L} w(x_j,x_i) v(x_i)=\theta v(x_j).
\end{equation}
If $L$ is small enough such that this eigenvalue problem can be solved, then $\theta$ and $v(x_i)$, $i=1,\ldots, L$, can be computed, which in turn can be substituted back into \eqref{nystrom2} to find an approximation to $v(y)$, for any $y\in V$. In short, we approximate the eigenvectors in \eqref{eig:eq} by extensions of the eigenvectors in \eqref{nystrom3}, using the extension equation \eqref{nystrom2}, and we approximate the eigenvalues in \eqref{eig:eq} by the eigenvalues from \eqref{nystrom3}. The main Nystr\"om assumption is that these approximated eigenvectors and eigenvalues approximately diagonalise $W$.
\red{For further details on the Nystr\"{o}m method, we refer the reader to Appendix \ref{nystrom:app} where a description of the method is given in matrix notation.}

\subsection{Pseudocode}  \label{subsec:algGL}

We present here the pseudocode combining all the different steps described above for the realisation of the GL minimisation. We recall that $\varepsilon$ is the scale parameter of the GL functional \eqref{def:GL_discr}, $\sigma$ is the variance used in the Gaussian similarity function \eqref{similarity_func}, $C$ is the convex splitting parameter in \eqref{convexsplit1}-\eqref{convexsplit2} and $L$ is the number of sample points in \eqref{nystrom1}.

\begin{algorithm}[h!]
\caption{GL-minimisation with Nystr{\"o}m extension for image segmentation}
\label{alg:GLmin_nystrom}
\begin{algorithmic}[1]
\STATE{\ \sc{Parameters}: \normalfont{$L\ll S$, $\sigma$, $\varepsilon$,  $C$}}.
\STATE {\ select $L$ random points and build the set $X\subset V$}
\STATE{\ get a partition $V=X\cup Y,~ Y:=X^c$}
\STATE {\ determine features and edge weights of $X$ and $Y$ using \eqref{similarity_func} and build $W_{XX}$ and $W_{XY}$} 
\STATE{\ \red{\textbf{Nystr{\"o}m extension} to compute normalised matrix of eigenvectors of $W$ and get eigenvalues-eigenvectors of $W$ $(\hat{\lambda}_i,v_i)$}} 
\STATE{\ output $\leftarrow$ eigenvalues-eigenvectors $(1-\hat{\lambda}_i,v_i)$ of $L_s$ used as GL minimisation input}
\STATE{\ \textbf{convex splitting} for GL minimisation through Fourier transform methods, as described in Section \ref{subsec:convexsplit}}
\STATE{\ output $\leftarrow$ the binary segmentation.}
\end{algorithmic}
\end{algorithm}

We will now give further details. First we randomly select $L$ pixels from $I$. As described in Section~\ref{subsec:graph} we now create a vertex set $V\cong I$, which we partition into a set $X$, consisting of the vertices corresponding to the $L$ randomly chosen pixels, and a set $Y:= V \setminus X$. We now compute the feature vectors of each vertex in $V$. If $I$ is a grey scale image, we can represent features by an intensity map $f:V \to \R$. If $I$ is an RGB colour image instead, we use a vector-valued (red, green, and blue) intensity map $f:V\to\R^3$ of the form $f(x)=(f_{\textcolor{red}{R}}(x), f_{\textcolor{green}{G}}(x), f_{\textcolor{blue}{B}}(x))$. We mirror the boundary to define neighbourhoods also on the image edges. The feature function $\psi: V \to \R^K$ concatenates the intensity values in the neighbourhood $\nu(x)$ of a pixel into a vector: $\psi(x) := (f(\nu_1(x)), \ldots, f(\nu_{\tilde \tau}(x)))^T$, where $\nu(x)=(\nu_1(x),\ldots,\nu_{\tilde\tau}(x))\in \R^{\tilde{\tau}}$ is the neighbourhood vector of $x\in V$ \red{defined in Section \ref{subsec:graph}} and $\tilde \tau = (2\tau+1)^2$, the size of the neighbourhood of $x$. Note that $K=\tilde \tau$ if $I$ is a grey scale image and $K=3\tilde \tau$ if $I$ is an RGB colour image.


Additional features can be considered, such as texture, for instance. For instance, we consider the eight MR8 filter responses \cite{varma} as texture features on a grey scale image and choose the function $t:V\to \R^8$ as $t(x) = (\text{MR8}_1(x), \ldots, \text{MR8}_8(x))$. Hence, the feature function $\psi$ is now defined as $\psi(x):=(t(\nu_1(x)),\ldots,t(\nu_{\tilde\tau}(x)))^T$ where $\nu(x)$ and $\tilde\tau$ are defined as above. Here, $K=8\tilde\tau$. Of course, a combination of colour and texture features can be considered as well by considering $\psi$ defined as $\psi(x):=(f(\nu_1(x)),t(\nu_1(x)),\ldots,f(\nu_{\tilde\tau}(x)),t(\nu_{\tilde\tau}(x)))$ for every $x$ in $V$. In this case, when dealing with RGB colour images, the dimension of the feature vector is therefore $K=11\tilde\tau$.


Using \eqref{similarity_func}, 
the Nystr\"om extension can be performed for approximating the eigenvectors and eigenvalues of
$W$ as described in Section \ref{subsec:nystrom} and in Appendix \ref{nystrom:app},
which are then used to compute the eigenvectors $\left\{ v_k \right\}$ of $L_s$ and corresponding eigenvalues $\left\{ \lambda_k\right\}$, compare \eqref{eigenvectorsLs}. 
Recalling \eqref{spectrdecun}, those eigenvectors are used as basis functions for $U_n$, the numerical approximation of $u$ in the $n$-th iteration of the GL minimisation. Considering \eqref{GLconvexsplit} and writing the nonlinear quantities appearing in terms of $\left\{ v_k \right\}$ similarly as in \eqref{spectrdecun}, we have for $x\in V$
\begin{align*}
&\left(U_n(x)\right)^3=\sum_{k} \beta_n^k(x)v_k(x),\quad\chi(x)\left(U^n(x)-u_0(x)\right)=\sum_{k} \gamma_n^k(x)~ v_k(x).
\end{align*}
The computation of $U$ in the next iteration reduces to finding the coefficients $\alpha_{n+1}^k$ in the expression
\begin{equation*}
 U_{n+1}(x)=\sum_{k} \alpha_{n+1}^k(x) v_k(x),\quad x\in V,
\end{equation*}
in terms of $\beta_n^k, \gamma_n^k$ and the other parameters involved, that is the scale parameter $\varepsilon$ in \eqref{def:GL_discr}, the parameter $C>0$ appearing in the splitting \eqref{GLdconvexsplit} and the time step $\Delta t$.
Using \eqref{GLconvexsplit}, we compute $ \alpha_{n+1}^k$ simply as
\begin{equation*}
 \alpha_{n+1}^k=\mathcal{D}_k^{-1}\left( \left(1+\frac{\Delta t}{\varepsilon} + C\Delta t\right)  \alpha_{n}^k - \frac{\Delta t}{\varepsilon}  \beta_{n}^k - \Delta t~ \gamma_{n}^k \right),
\end{equation*}
where $\mathcal{D}_k$ is defined as $\mathcal{D}_k:= 1 + \Delta t(\varepsilon\lambda_k + C)$. 
 

\section{Hough transform for scale detection}    \label{sec:hough_use}

In order to detect objects in an image with specific, \emph{a priori} specified shapes, in the following we will make use of the Hough transform. For our purposes, we will focus in particular on straight lines detection (for which the Hough transform was originally introduced and considered \cite{hough}) and circles, \cite{dudahart}. Other applications of this transformation for more general curves exist as well. In \cite{beltrampiana,campi1} the Hough transform is used in the context of astronomical and medical images for a specific class of curves (\emph{Lamet} and \emph{Watt} curves). In \cite{grah} applications to cellular mitosis are presented. There, the Hough transform recognises the cells (as circular/elliptical objects) and tracks them in the process of cellular division. 
\red{For more details on the use of the Hough transform for line and circle detection we refer the interested reader to Appendix \ref{sec:hough}.}

\paragraph{Numerical strategy.} \red{Hough transform methods for edge detection are usually applied to binary images. Therefore, we start by using the classical Canny method for edge detection \cite{canny1} in which we replace the original preliminary Gaussian filtering by an edge-preserving Total Variation smoothing \cite{ROF} which has the advantage of removing noise while preserving edges. This step will result in a binary image for the most prominent edges in the image}. Having decided which geometrical shape we are interested in (and, as such, its general parametric representation), the corresponding parameter space is subdivided into accumulator arrays (cells) whose dimension depends on the dimension of the parameter space itself (2D in the case of straight lines, 3D in the case of circles). Each accumulator array  groups a range of parameter values. The accumulator array is initialised to $0$ and incremented  every time an object in the parameter space passes through the cell. In this way, one looks for the peaks over the set of accumulator arrays as they indicate a high value of intersecting objects for a specific cell. In other words, they are indicators of potential objects having the specific geometrical shape we are interested in.

\subsection{Pseudocode} \label{subsec:hough}

Numerically, dealing with the Hough transform consists of looking for peaks of the accumulator arrays in the parameter space onto which the original image is mapped. We use the MATLAB routines \texttt{hough}, \texttt{houghpeaks}, and \texttt{houghlines} for straight lines detection and \texttt{imfindcircles} for circle detection. The accuracy and the number of detections for such routines can be tuned by some parameters, such as, for instance, the maximum number of peaks one wants to consider, $obj_{max}$, or the array peak threshold, $thresh$, i.e. the minimum number of elements for an accumulator array to be considered a peak. The user also has to specify an initial range of pixel values $[s_{min},s_{max}]$ as a very rough approximation of the measurement scale. Namely, in the case of line detection this determines a minimum/maximum \emph{spacing} between lines , whereas for circle detection this serves as a rough approximation of the range of values for the circles' radii.  This rough approximation may be given, for example, from average data which the user knows \emph{a priori}. We explain this with some examples in Section \ref{sec:res}. Accuracy of the detection algorithm is tuned by a parameter $acc$. In case of linear objects this corresponds to choose the maximum number of pixels between two line segments to consider them as one single line, whereas for circle detection this corresponds to the circularity of an object to be considered a circle.

\begin{algorithm}
\caption{Hough transform for lines and circles detection}
\label{alg:Hough}
\begin{algorithmic}[1]
\STATE{\ \sc{Parameters}: \normalfont{ $[s_{min}, s_{max}]$, $obj_{max}$, $acc$, $thresh$} }
\STATE {\ preprocessing: \red{TV-Canny} edge detection}
\STATE {\ compute the Hough transform of the edge image}
\STATE{\ set up detection accuracy, depending on $acc$, and use $[s_{min}, s_{max}]$ as rough initial guess}
\STATE{\  determine at most $obj_{max}$ peaks in the parameter space, thresholding using $thresh$ }
\STATE{\ output $\leftarrow$ peaks in the parameter space, corresponding to objects of interest in the original image}
\end{algorithmic}
\end{algorithm}

\section{Method, numerical results, and applications} \label{sec:res}

We report in this section the numerical results obtained by the combination of the methods presented for the detection and quantitative measurement of objects in an image.

To avoid confusion, we will distinguish in the following between two different meanings of  scale. Namely, by image scale we denote the proportion between the real dimensions (length, width) of objects in the image and their corresponding dimensions quantified in pixel count. Dealing with measurement tools, we talk about measurement  scale to intend the ratio between a fixed unit of measure ($mm$ or $cm$) marked the measurement tool considered  and the correspondent number of pixels on the image. 

\subsection{Male pied flycatcher's blaze segmentation} \label{sec:flycatcher}

Here we present the numerical results obtained by applying algorithms \ref{alg:GLmin_nystrom} and \ref{alg:Hough} to the male pied flycatcher blaze segmentation and measurement problem described in the introduction. Our image database consists of 32 images of individuals from a particular population of flycatchers. Images are $3648 \times 2736$ pixels and have been taken by a Canon 350D camera with Canon zoom lens EFD $18$--$55\, mm$, see Figure \ref{fig:im_database}. In each image one of two types of measurement tool is present: a standard ruler or a surface on which two concentric circles of fixed diameter ($1\, cm$ the inner one, $3\, cm$ the outer one) are drawn. In the following we will refer to these tools as \emph{linear} and \emph{circular} ruler, respectively. Here,  the measurement scale corresponds to the distance between ruler marks for linear rulers and to the radius of the inner circles for circular rulers.

Figure \ref{fig:im_database} shows clearly that the scale of the images in the database may vary significantly because of the different positioning of the camera in front of the flycatcher. In order to study possible correlations between the dimensions (i.e. perimeter, area) of the blaze and significant behavioural factors, the task then is to segment the blaze and detect automatically the scale of the image considered to provide scale-independent measurements. 

\red{\paragraph{Parameter choice for Algorithm \ref{alg:GLmin_nystrom}} The GL-segmentation method exploits similarities and differences between pixels in terms of RGB intensities and texture within their neighbourhood. In our image database these similarities and differences are very distinctive and will guide the segmentation step. 
Recalling Section \ref{subsec:algGL}, we note that some parameters need to be tuned for the graph GL minimisation. Those are the number $L$ of Nystr\"{o}m points, the variance $\sigma$ of the similarity function \eqref{symmgraphlapl}, the GL parameter $\varepsilon$ and the parameter $C$ for the convex splitting \eqref{GLdconvexsplit}. However, in our numerical experiments we had to tune these parameters only once. Namely, regarding the choice of $L$ for both the head and blaze segmentation, we used values not bigger than $5\%$ of the total size of the image considered. The variance appearing in the similarity function \eqref{similarity_func} was set to $\sigma^2=20$ and the weighting parameter $\varepsilon$ was chosen as $\varepsilon=0.01$ (a smaller choice would create numerical instabilities) and we set the convexity parameter $C=25$ or larger in order to guarantee the convexity of the functional appearing in \eqref{convexsplit2}.}

\red{\paragraph{Parameter choice for Algorithm \ref{alg:Hough}}: We briefly comment also on the choice of the parameters for the Hough transform, that is Algorithm \ref{alg:Hough}. Depending on the type of measurement tool considered (linear or circular ruler), different parameter selection methods are considered. In the case of linear rulers: for the longest line detection (i.e. ruler edge identification) the parameters $obj_{max}$ and $thresh$ were set to $1$ and to $85\%$ of the maximum value of the Hough transform matrix, respectively; for the detection of the ruler notches, the same parameters were chosen as  $obj_{max}=500$ and $thresh=20\%$ of the maximum value of the Hough transform matrix were used. As discussed in \ref{subsec:hough}, the range $[s_{min}, s_{max}]$  was chosen based on previously collected average data on the head diameter. In particular, after the head detection step, the number of pixels corresponding to the diameter of the head was automatically computed by means of the option \texttt{EquivDiam} of the MATLAB routine \texttt{regionprops} and compared with the average measurement of $1.51$ $cm$ provided by available databases on pied flycatcher. In this way, an initial, rough estimate of the ruler scale is found and used to determine a \emph{spacing parameter} $s$ and the interval $[s_{min}, s_{max}]$ by setting $s_{min}=s/2$ and $s_{max}=2s$. This range serves as a \emph{suppression} neighbourhood: once a peak in the Hough transform matrix (i.e. a line or a circle) is identified, starting from it the successive peaks found outside this range are set to $0$ (i.e. possible lines or circles within this interval are discarded), while only the ones inside the range (typically, the following line/circle we want to detect) are kept. For our problem, this corresponds to identifying as candidates for ruler notches only the lines away from each other at least $s_{min}$ and at most $s_{max}$ from the peak which has been previously identified. Analogously, the same can be done with circular rulers, where we recall the inner/outer radii are $1$ $cm$ and $3$ $cm$, respectively. In this case $obj_{max}=2$ since the circular ruler is made of only two concentric circles.}

\paragraph{Comparison with Chan-Vese model.} Due to the very irregular contours and the fine texture on the flycatcher's forehead, standard variational segmentation methods such as Canny edge detection or Mumford-Shah models, \cite{mumfshah,ambrosiotort,braides,braidesdalmaso},  are not suitable for our task, as preliminary tests showed. Chan-Vese \cite{chanvese}
 is not suitable either, mainly because of the small scale detection limits, the dependence on the initial condition, and the parameter sensitivity which may prevent us from an automatic and accurate segmentation of the tiny, yet characteristic feathers composing the blaze. In particular, the optimal parameters $\mu$ and $\nu$ appearing in the Chan-Vese functional and a sufficiently accurate initial condition need to be chosen typically by trial and error for every image at hand.

\begin{figure}[!h]
\begin{center}
\begin{subfigure}{.23\textwidth}
\centering
\includegraphics[height=3cm]{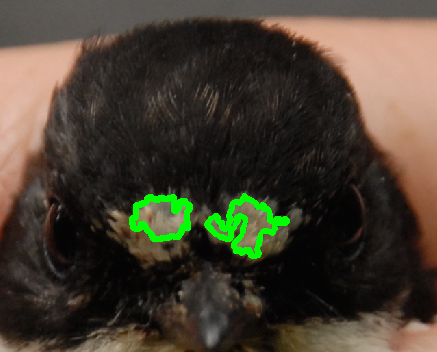}
\caption{Chan-Vese segmentation}
\label{fig:skinmoles_detection1}
\end{subfigure}
\begin{subfigure}{.23\textwidth}
\centering
\includegraphics[height=3cm]{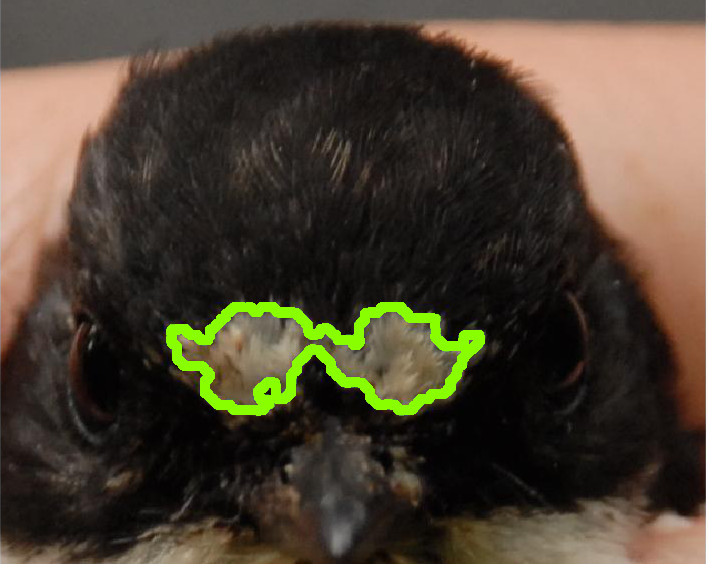}
\caption{GL-segmentation}
\label{fig:skinmoles_detection2}
\end{subfigure}
\caption{Blaze segmentation results computed by using Chan-Vese model \cite{chanvese} and  GL minimisation  (Algorithm \ref{alg:GLmin_nystrom}). The dependence of the Chan-Vese model on the initial condition and its sensitivity to the model parameters may result in inaccurate detections, while the GL approach provides more reliable segmentation results.}
\label{fig:comparison_chan_vese}
\end{center}
\end{figure}

For comparison, we report in  Figure \ref{fig:comparison_chan_vese} the blaze segmentation results obtained by using Chan-Vese model (see \cite{chanvese,chansandvese}) and our graph based method which will be described in more detail in the following.

\subsubsection{Detailed description of the method} 
\begin{figure}[!h]
\begin{center}
\includegraphics[height=3.2cm]{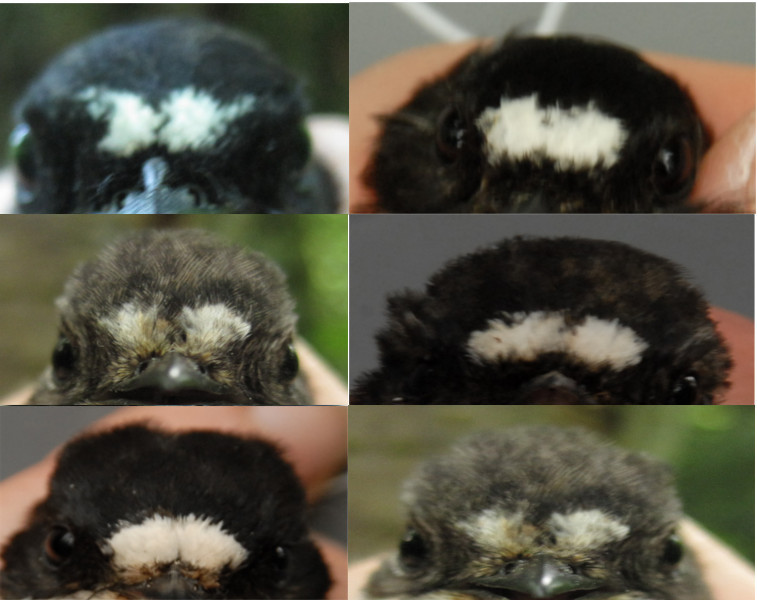}
\includegraphics[height=3.2cm]{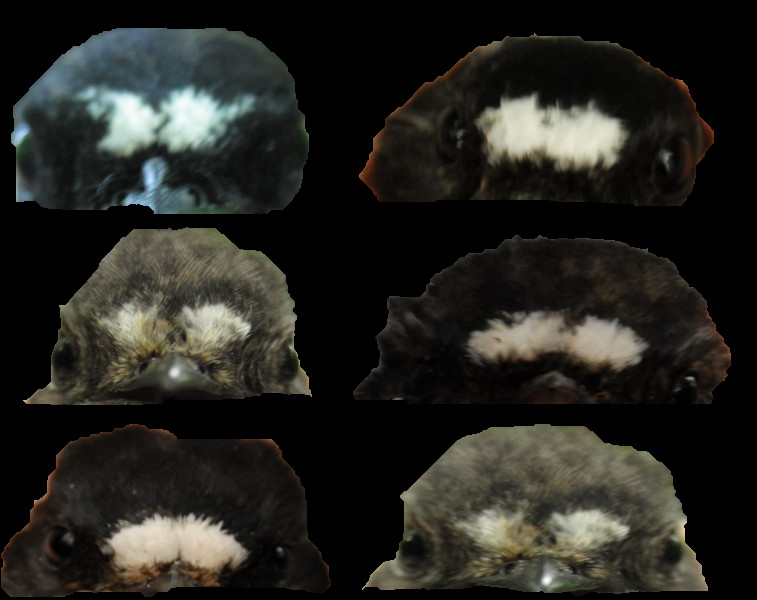}
\caption{Training dictionary for head detection: the heads are manually selected by the user and separated from the background. Then, the corresponding regions are labelled with $1$ while the background is labelled by $-1$.}
\label{fig:dictionary}
\end{center}
\end{figure}

We divide our task into different steps:
\begin{enumerate}
\item For a given, unsegmented image, we detect the head of the pied flycatcher through a comparison with a user prepared dictionary  (see Figure \ref{fig:dictionary}) using GL segmentation Algorithm \ref{alg:GLmin_nystrom}. Further computations are restricted to the head only.
\item  Starting from the reduced image, a second step similar to Step 1 is now performed for the segmentation of the blaze, using again Algorithm \ref{alg:GLmin_nystrom}. A  dictionary of blazes is used an extended set of features is considered.
\item A refinement step is now performed in order to reduce the outliers detected in the process of segmentation. 
\item We use the Hough transform based Algorithm~\ref{alg:Hough} to detect in the image objects with \emph{a priori} known geometrical shape (lines for linear rulers, circles for circular rulers) for the computation of the measurement scale. 
\item The final step is the measurement of the values we are interested in (i.e. the perimeter of the blaze, its area and the width and height of the different blaze components).  In the case of linear rulers our results are given up to some error (due to the uncertainty in the detection of the measurement scale computed as average between ruler marks distances).
\end{enumerate}

\red{Figure \ref{fig:diagram} shows a diagram which outlines the workflow of our method.
\begin{figure*}[!h]
\begin{center}
\includegraphics[width=\textwidth]{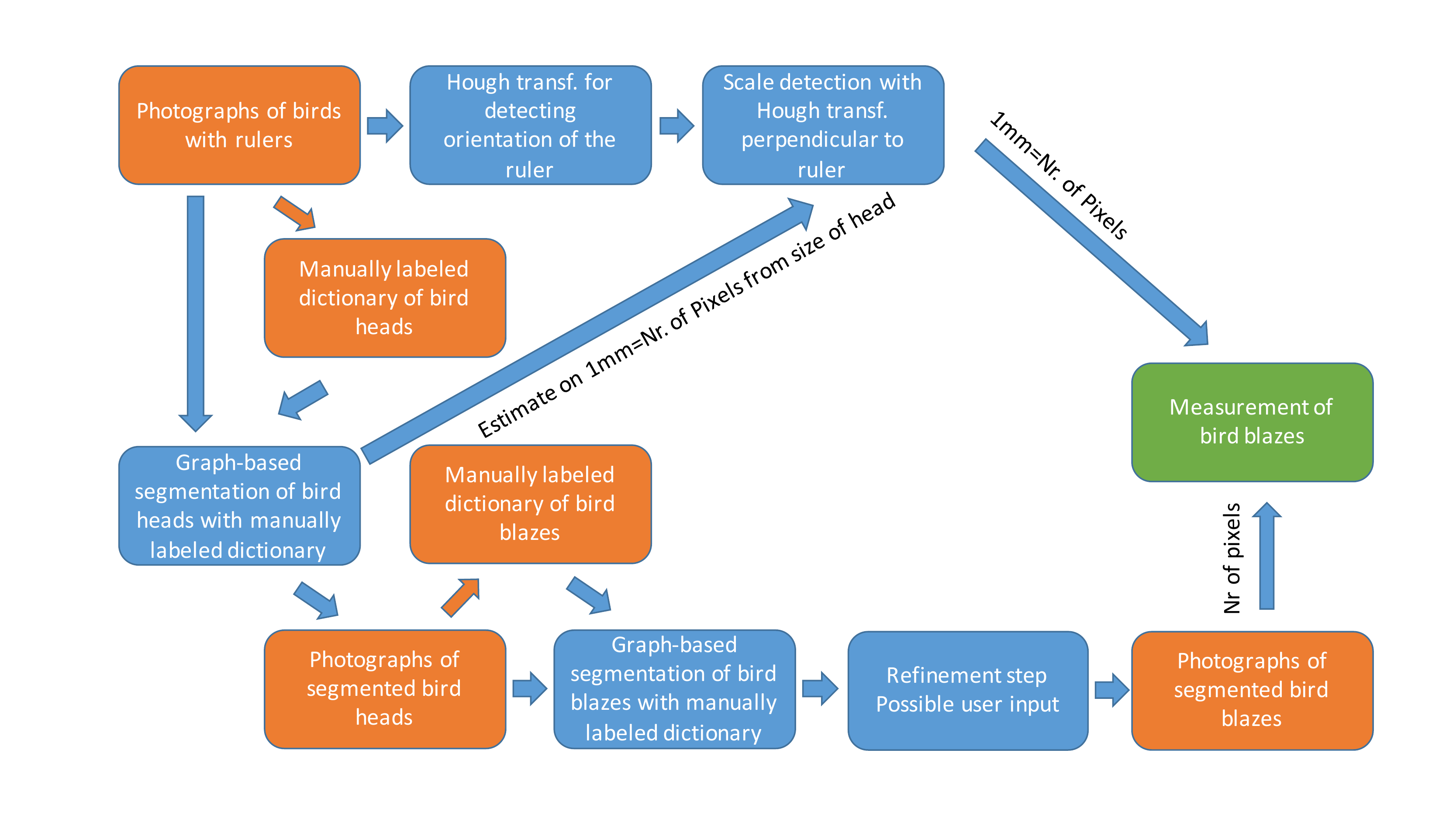}
\end{center}
\caption{The diagram describes the different steps of the segmentation/measurement procedure. Boxes requiring the user input are coloured orange, while the ones where the automatic segmentation/measurement steps are performed are coloured blue. The final objective is coloured green.}
\label{fig:diagram}
\end{figure*}
}
In order to establish relations with behavioural and biological data confirming or contradicting the initial assumption of correlation between blaze size and higher attractiveness  presented in the introduction \cite{pottimontalvo}, we have implemented a user ready program for the quantitative analysis and measurements of the size of the bird blazes which is currently used by the Department of Zoology of the University of Cambridge. The results of this study will be the topic of a forthcoming paper \cite{ourpreprint}.

\medskip

In the following we give more details about each step.

\paragraph{Step 1: Head detection.} We consider unlabelled images in the database and compare each of them with a dictionary of previously labelled images, see Figure \ref{fig:dictionary}. The training regions (i.e. the heads) are labelled with a value $1$, the background with value $-1$. Unlabelled regions are initialised with value $0$.

The main computational difficulties in this step are due to size of the images considered. This may affect the performance of the algorithm as in order to apply the Nystr{\"o}m completion technique described in Section \ref{subsec:nystrom} one has to choose an adequate number of points whose features will approximate well the whole matrix. The larger and more heterogeneous the image is, the larger will be the number of points needed to produce a sensible approximation. We circumvent this issue noticing that at this stage of the algorithm, we only need a rough detection of the head which will be used in the following for the accurate segmentation step. Thus, downscaling the image to a lower resolution (in our practice, reducing the resolution by ten times the original one) allows us to use a small number of Nystr{\"o}m sample points (typically $150$--$200$) to produce an accurate result.

\begin{figure}[!h]
\begin{center}
\includegraphics[height=3.1cm, width=3.8cm]{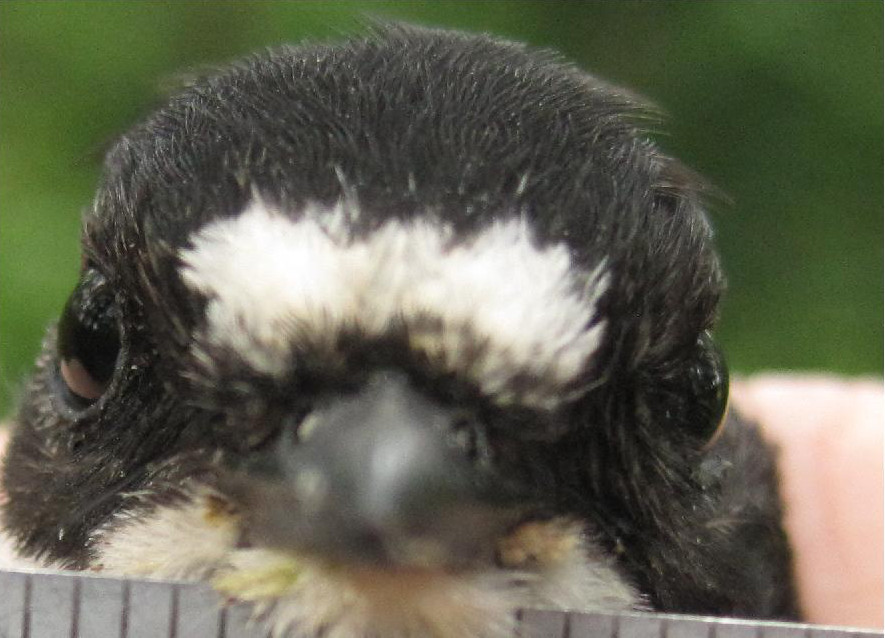}\hspace{0.2cm}
\includegraphics[height=3.1cm, width=3.6cm]{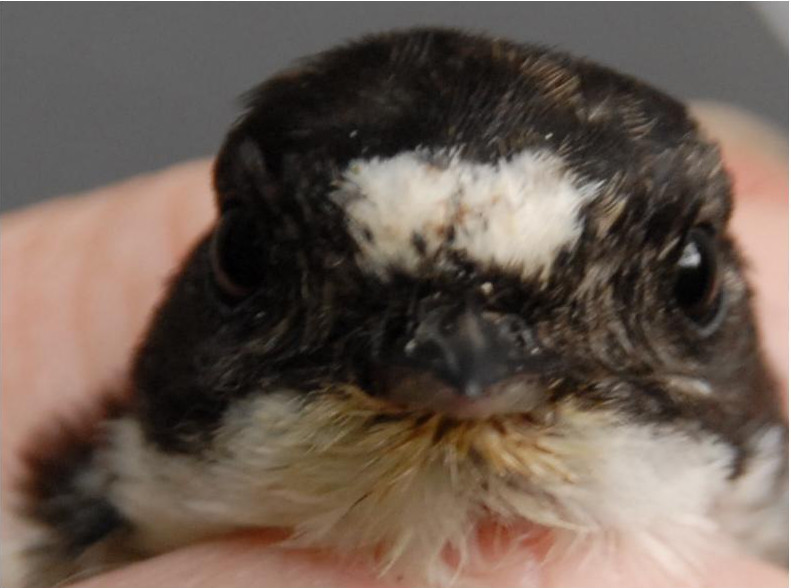}
\caption{Head detection from images in  Figure \ref{fig:im_database} using dictionary
in Figure \ref{fig:dictionary}}
\label{fig:detected_head}
\end{center}
\end{figure}

For this first step we use as features simply the RGB intensities and proceed as described in Section~\ref{subsec:algGL}. Once the head is detected, the resulting image is upscaled again to its original resolution. The solutions computed for the images in Figure \ref{fig:im_database} are presented in Figure \ref{fig:detected_head}. 

\begin{figure}[!h]
\begin{center}
\includegraphics[height=3cm]{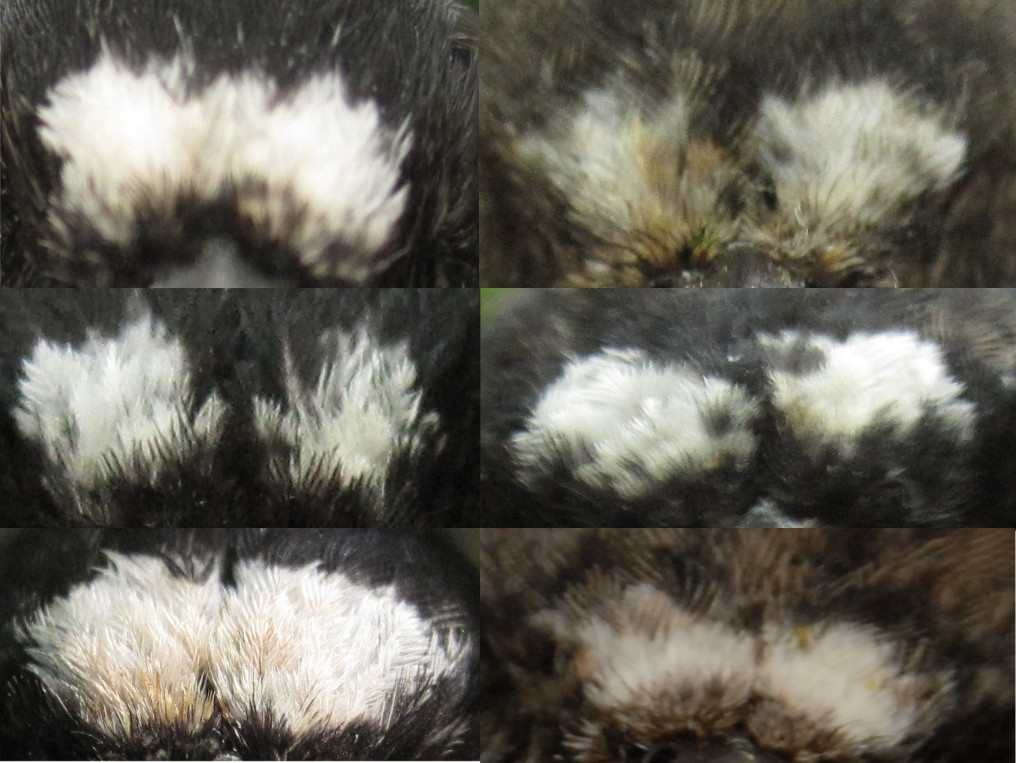}
\includegraphics[height=3cm]{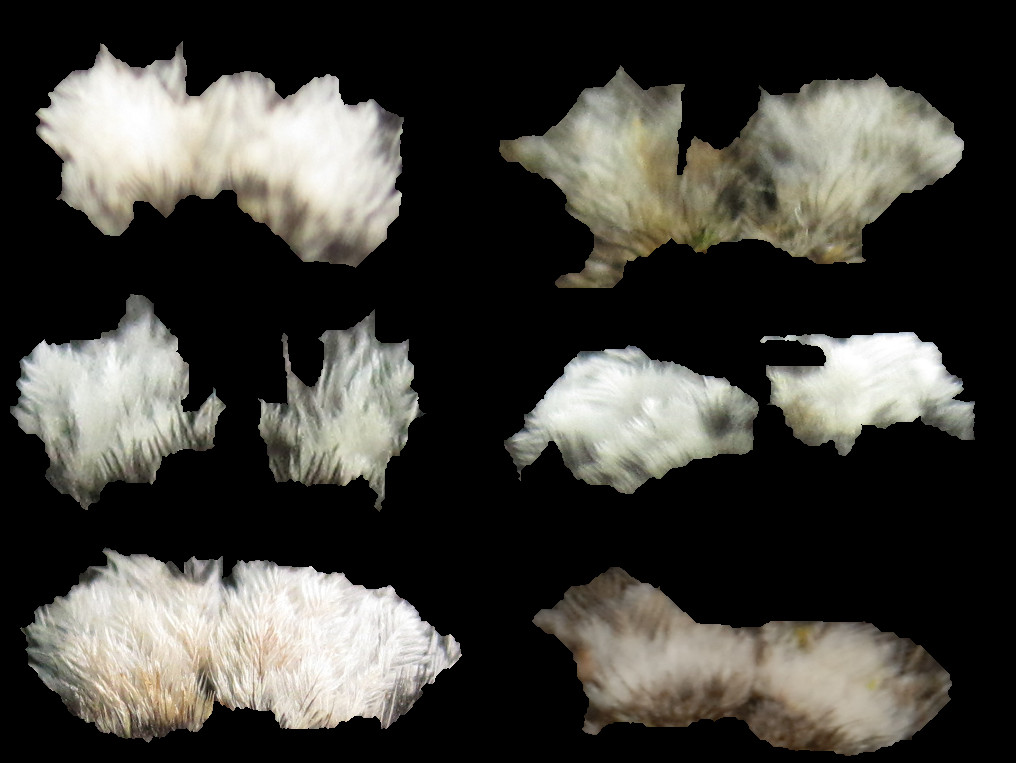}
\end{center}
\caption{Training dictionary for blaze segmentation. As in Figure \ref{fig:detected_head} blazes are manually selected by the user and labelled with $1$, while black feathers on the background are labelled with $-1$.}
\label{fig:dictionary_blazes}
\end{figure}

\begin{figure}[!h]
\begin{center}
\includegraphics[height=2.9cm]{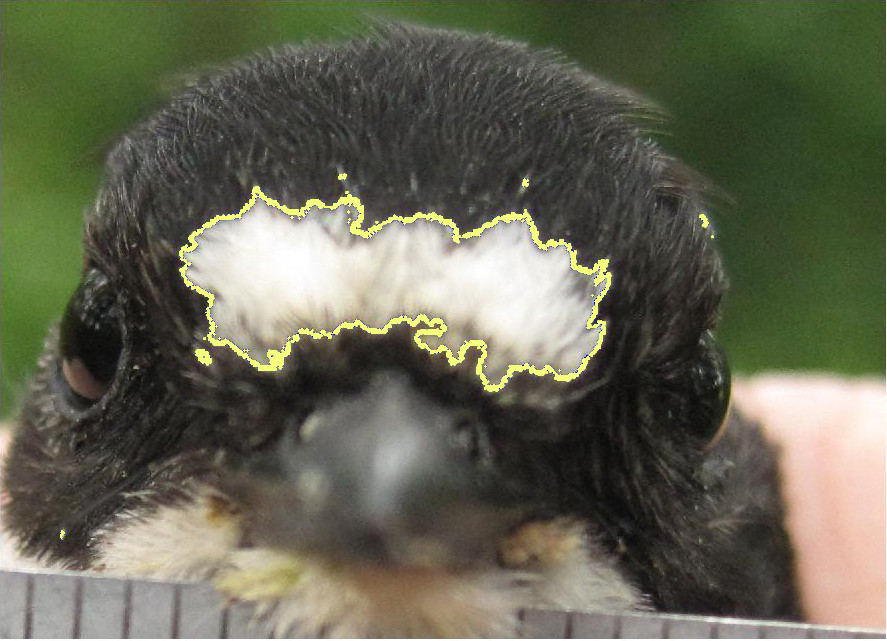}
\includegraphics[height=2.9cm]{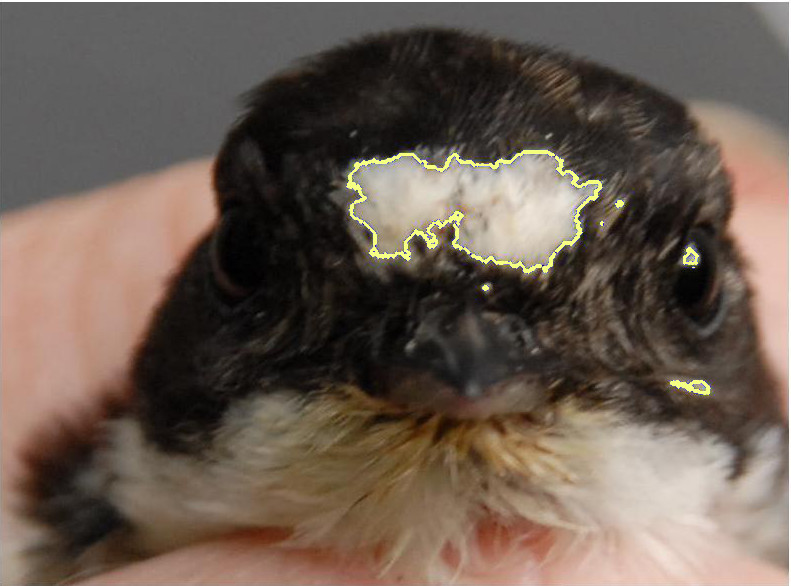}
\caption{Blaze segmentation}
\label{fig:blaze_segmentation}
\end{center}
\end{figure}

\paragraph{Step 2: Blaze segmentation.} We consider now the reduced image from which we want to extract the flycatcher's blaze. Again, a dictionary of different blazes is manually created by the user (see Figure \ref{fig:dictionary_blazes}). Again, training regions (the blazes) are labelled with value $1$ and the black feathers in the background with value $-1$. As before, unlabelled regions are initialised with value $0$. At this stage, RGB intensities alone are not enough to differentiate the blazes from the background consistently in a large number of bird images, due to the colour difference between different blazes. For this step, an additional feature to be considered is the \textbf{texture} of the blaze.  For this purpose, we use the MR8 texture features presented in \cite{varma} and proceed as detailed in Section~\ref{subsec:algGL}. For $3\times3$ neighbourhoods, the feature vector for each pixel will be an element in $\R^{99}$, see Section \ref{subsec:algGL}. The Ginzburg-Landau minimisation provides the segmentation results shown in Figure \ref{fig:blaze_segmentation}.

\begin{figure}[!h]
\begin{subfigure}{.5\textwidth}
\centering
\includegraphics[width=6cm]{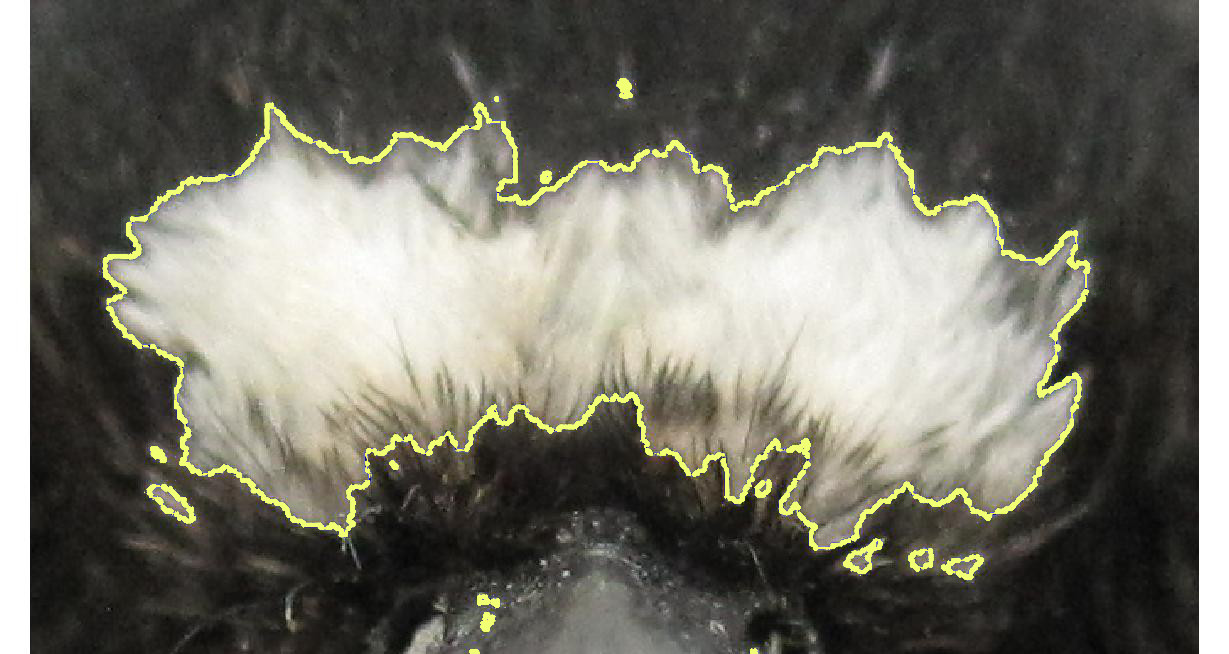}
\caption{Before refinement}
\label{fig:beforerefin}
\end{subfigure}
\begin{subfigure}{.5\textwidth}
\centering
\includegraphics[width=6cm]{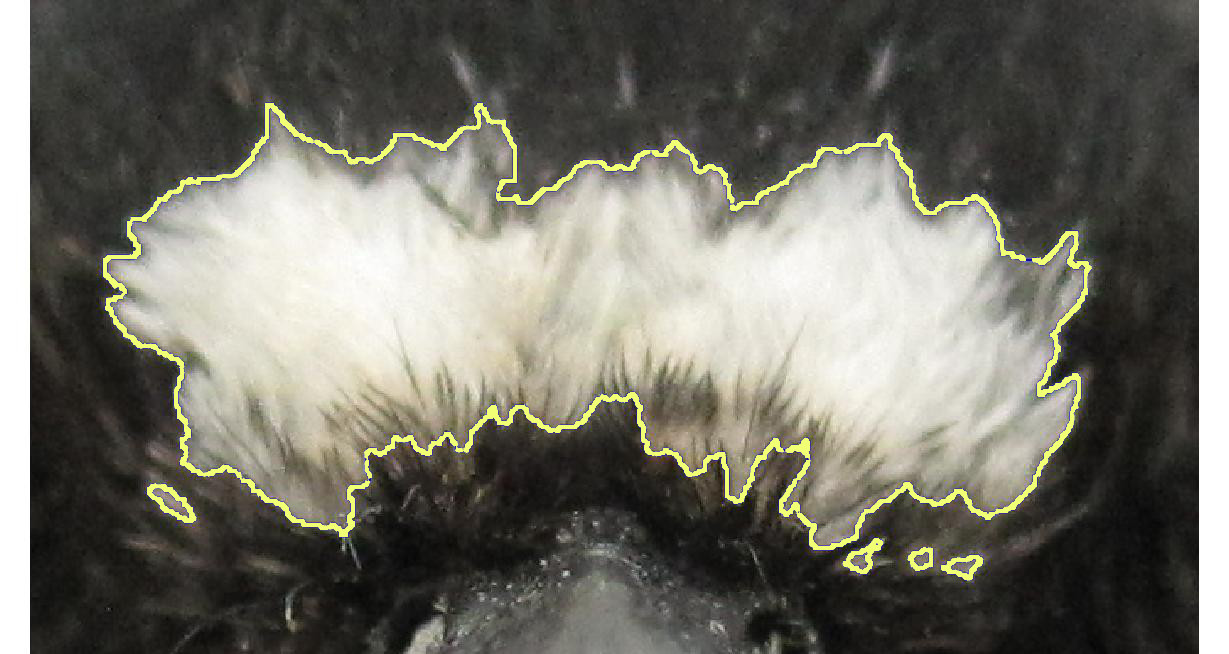}
\caption{After refinement}
\label{fig:afterrefin}
\end{subfigure}
\caption{Example of segmentation refinement}
\label{fig:refinementdiff}
\end{figure}

\begin{figure*}[!h]
\begin{center}
\includegraphics[width=15cm,height=4cm]{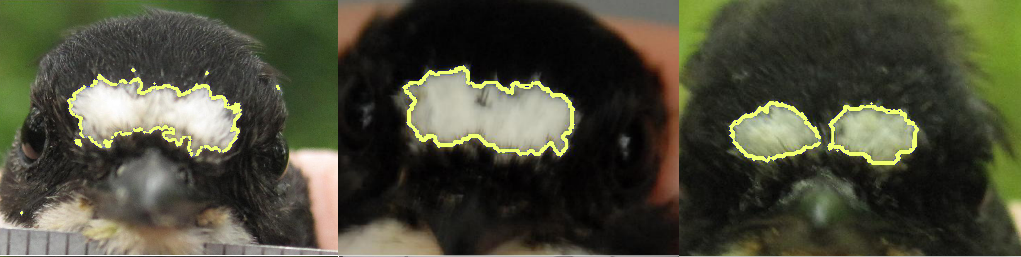}\vspace{0.3cm}
\includegraphics[width=15cm,height=4cm]{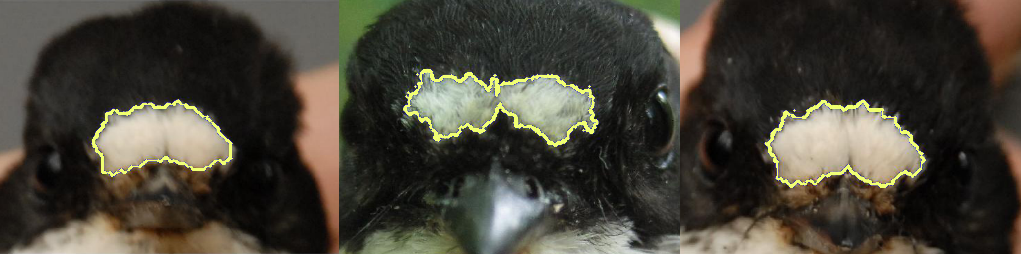}
\end{center}
\caption{Segmentation results after refinement step}
\label{fig:blaze_refinement}
\end{figure*}

\paragraph{Step 3: Segmentation refinement.}  This step uses very simple morphological operations in order to remove false detections obtained after Step 2. These can occur due to the choice of colour-texture based features  used to compute the feature vectors in Step 2. For instance, when looking at Figure \ref{fig:blaze_segmentation} (right) we observe that some bits on the left pied flycatcher's cheek have been detected as they exhibit similar texture properties as the ones on the blaze. In order to prevent this, our software asks the user to confirm whether the segmentation result provided is the expected one or if there are additional unwanted regions detected. If that is the case, using the MATLAB routine \texttt{bwconncomp} we label all the connected components segmented in the previous step, discarding among them all the ones whose area is smaller than a fixed percentage (we use $10\%$) of the largest detected component (presumably, the blaze). This works well in practice, see Figure \ref{fig:refinementdiff}. If the user is not satisfied he or she can remove manually the unwanted regions. Figure \ref{fig:blaze_refinement} shows some blaze segmentation results after the refinement step.

\red{\begin{remark}[Robustness to noise]
In order to reproduce the more realistic situation of images suffering from noise, we artificially added Gaussian noise with zero mean and different variances to some of the images in our database and performed the three analysis steps of our method. We report in Figure \ref{fig:noise_robustness} the results corresponding to two noise variances ($\sigma_1^2=0.02$, $\sigma_2^2=0.05$). The presence of noise influences both the head and blaze segmentation only slightly; the combination of RGB and texture features extracted in the neighbourhood of each point combined with the comparison to the dictionary make the algorithm  robust to noise and allows for an accurate blaze segmentation even in the noisy case.
\begin{figure}[!h]
\begin{center}
\begin{subfigure}{0.23\textwidth}
\includegraphics[height=2.9cm]{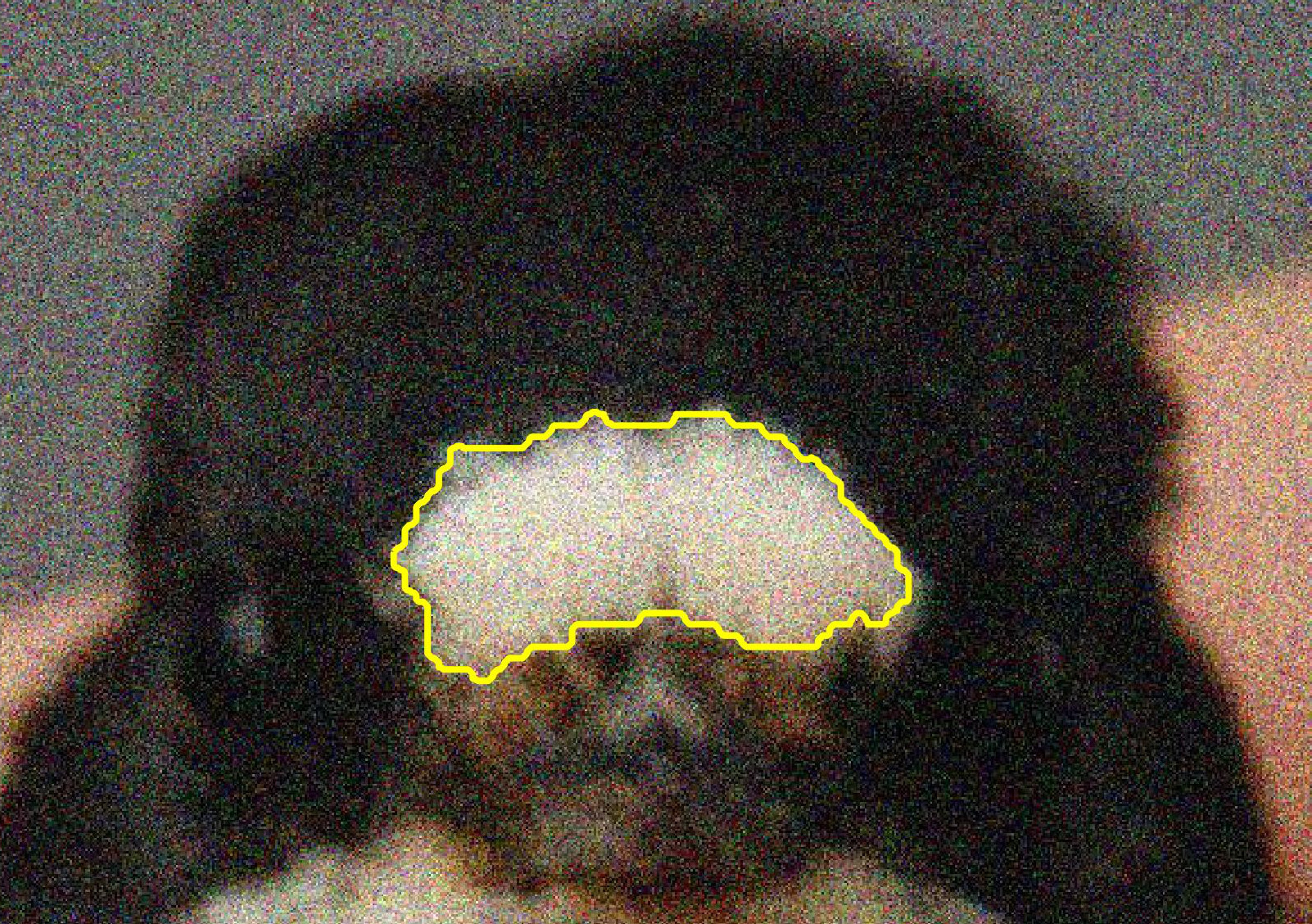}
\caption{$\sigma_1^2=0.02$}
\end{subfigure}
\hspace{0.1cm}
\begin{subfigure}{0.23\textwidth}
\includegraphics[height=2.9cm]{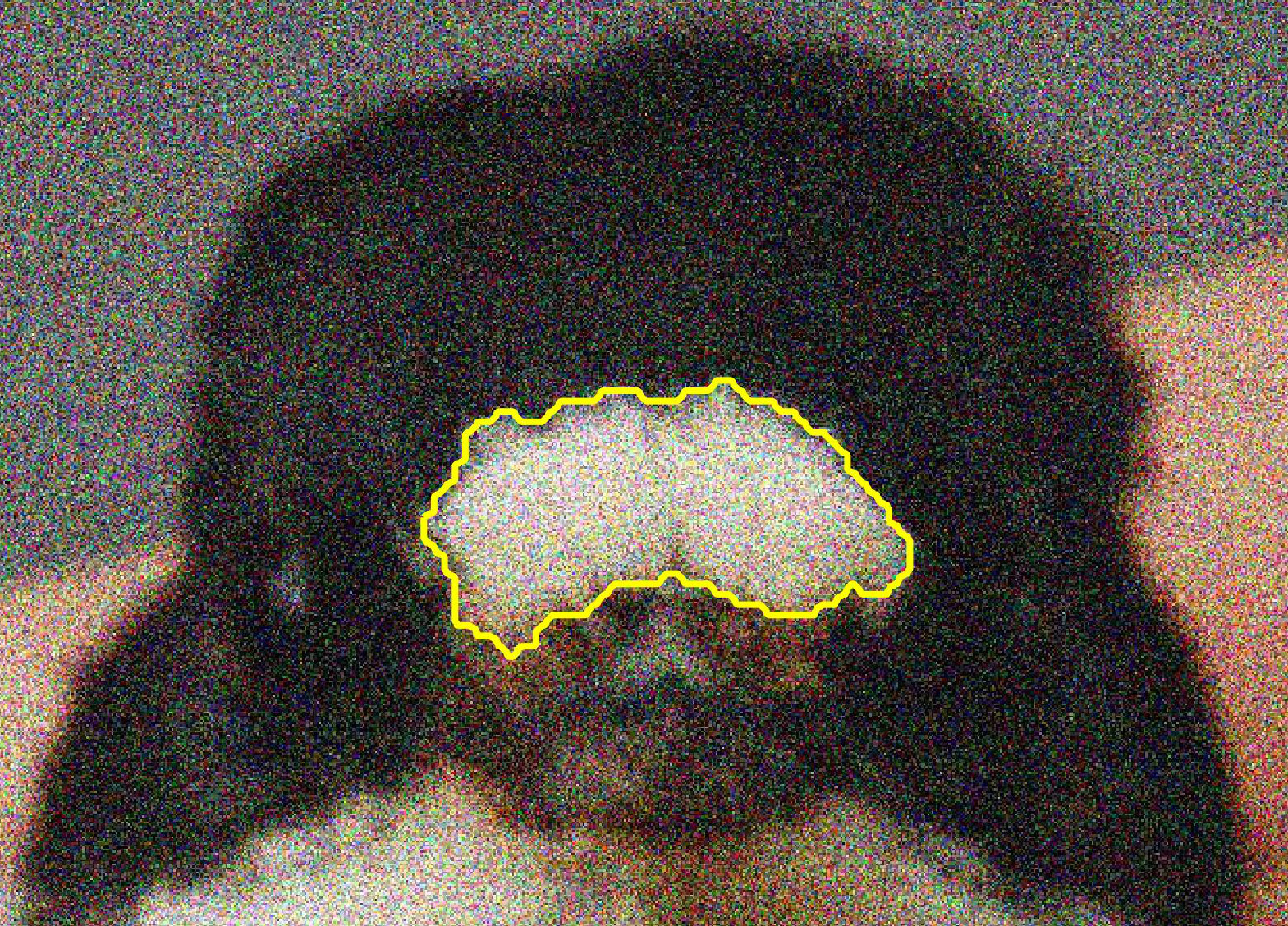}
\caption{$\sigma_1^2=0.05$}
\end{subfigure}
\caption{\red{Robustness to noise oscillations of GL minimisation for binary segmentation. Images have been artificially corrupted with Gaussian noise with zero mean and different variances.}}
 \label{fig:noise_robustness}
 \end{center}
\end{figure}
\end{remark}
\begin{remark}[Comparison with MBO segmentation]
We compare the blaze segmentation results obtained by minimising the discrete GL functional with the ones obtained using the segmentation algorithm considered in \cite{MerkujevKostic2013} as a variant of the classical Merriman-Bence-Osher (MBO) scheme \cite{MBO1992}. More details on the connections between this approach and the GL minimisation as well as some insights on its numerical realisation are given in Appendix \ref{appendix:MBO}. Following faithfully what is described in Section \ref{subsec:graph} and \ref{subsec:nystrom} for the graph and the operator construction step, respectively, we implemented the MBO segmentation algorithm following \cite[Section 2]{MerkujevKostic2013}. We remark that the MBO method has the advantage of eliminating the dependence on the interface parameter $\varepsilon$ of the GL functional by means of a combination of heat diffusion and a thresholding step. Instead of $\varepsilon$ the heat diffusion time $\tau$ needs to be chosen. In our numerical implementation we used $\tau=0.005$. Since no convex splitting strategies are required in this case, due to the absence of the non-convex double-well term, standard Fourier transform methods are used to solve the resulting time-stepping scheme. In Figure \ref{fig:MBO_comparison} we report the blaze segmentation results obtained after applying a refinement step similar to the one described above: we note that a segmentation result comparable to the ones shown in Figure \ref{fig:blaze_refinement} is obtained. Moreover, robustness to noise is observed also in this case. In terms of computational times, we observed that the replacement of the GL minimisation step with the MBO one did not affect significantly the speed of the segmentation algorithm.
\begin{figure}[!h]
\begin{center}
\begin{subfigure}{0.23\textwidth}
\includegraphics[height=2.9cm,width=3.9cm]{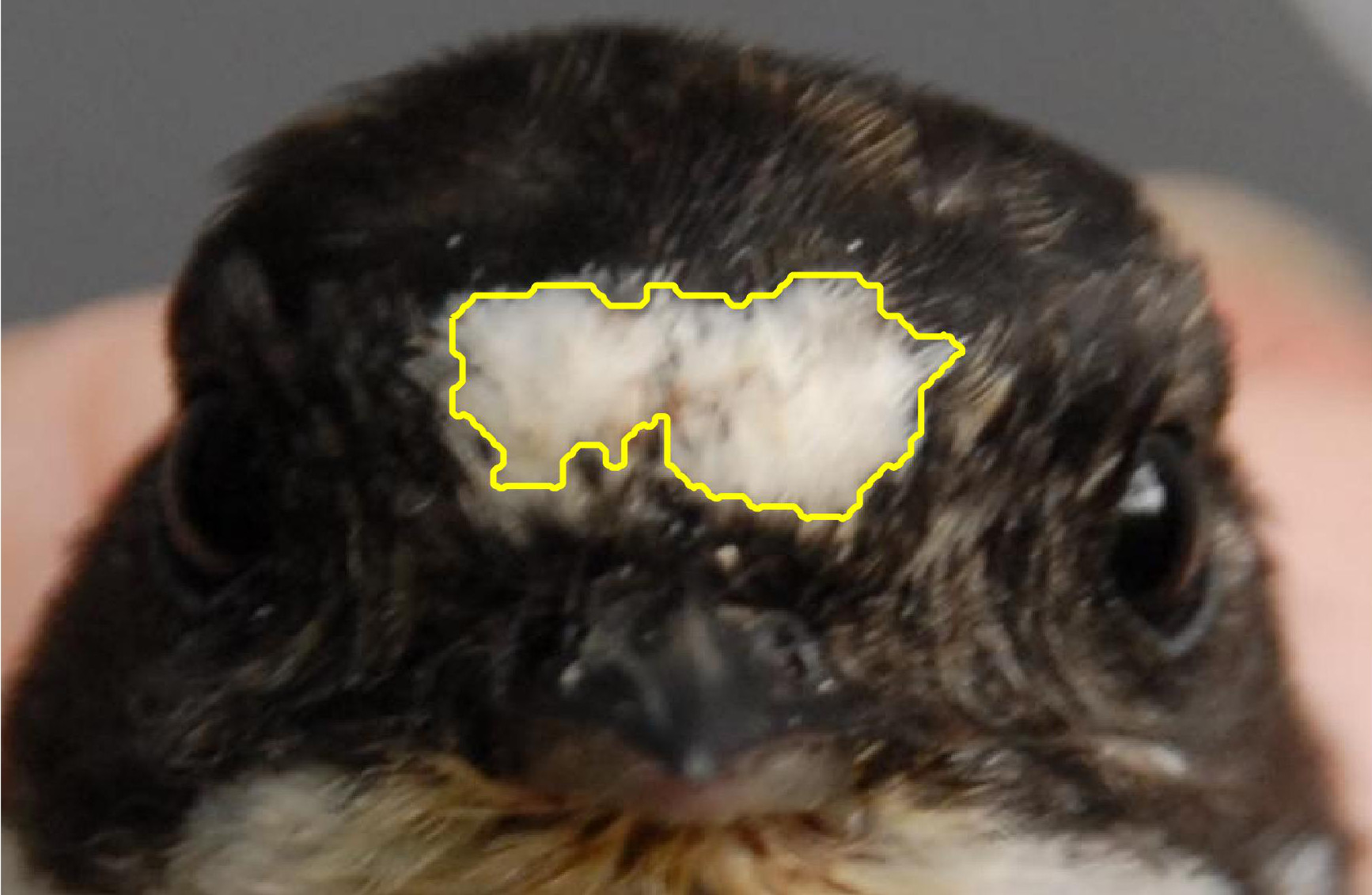}
\caption{MBO result}
\end{subfigure}
\begin{subfigure}{0.23\textwidth}
\includegraphics[height=2.9cm,width=3.9cm]{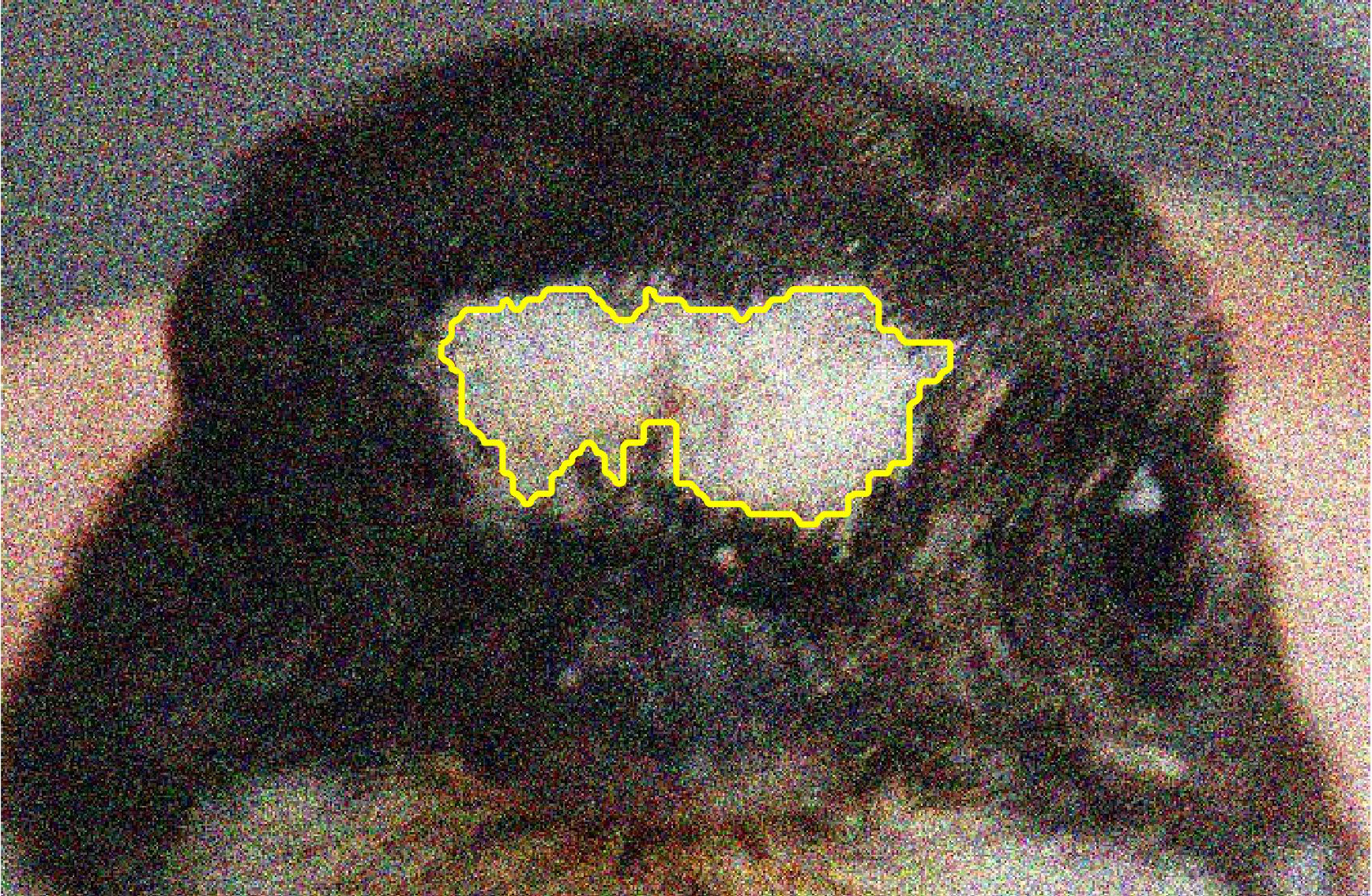}
\caption{MBO result, $\sigma^2=0.05$.}
\end{subfigure}
\caption{\red{Blaze segmentation results obtained by the MBO segmentation algorithm described in \cite{MerkujevKostic2013}, after refinement step. Robustness to noise is observed in this case as well. In both numerical tests, the diffusion time is chosen as $\tau=0.005$.}}
 \label{fig:MBO_comparison}
 \end{center}
\end{figure}
\end{remark}
}

\begin{figure}[!h]
\begin{center}
\includegraphics[height=2.5cm]{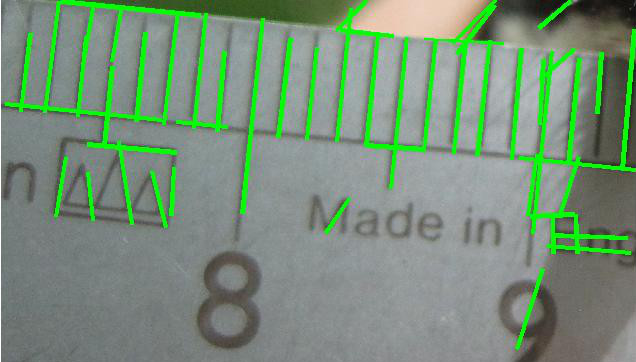} 
\includegraphics[height=2.5cm]{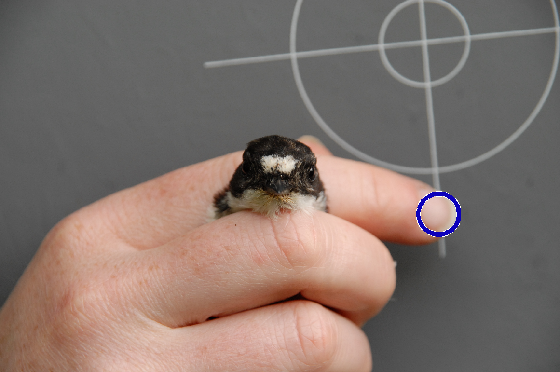}
\end{center}
\caption{Shadows, blur, noise or other objects in the image may disturb the detection. }
\label{fig:wrongdetections}
\end{figure}

\begin{figure}[!h]
\begin{center}
\includegraphics[height=3.1cm]{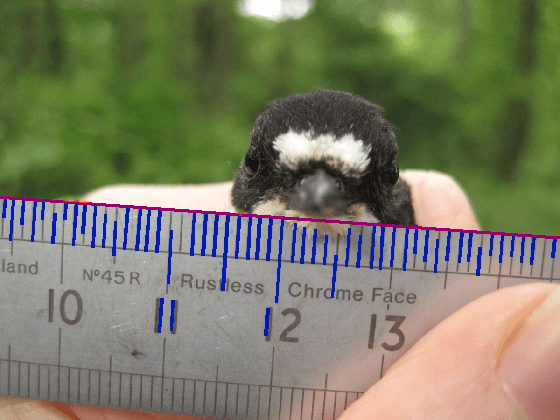}
\includegraphics[height=3.1cm]{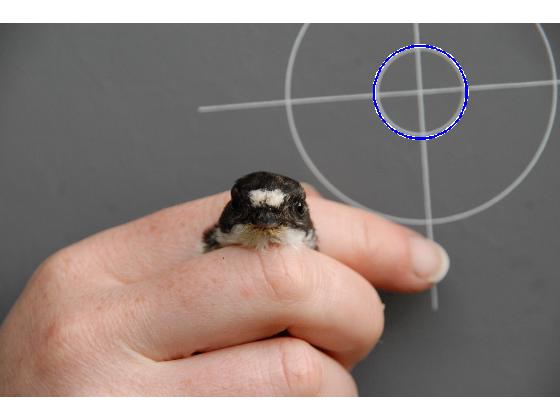}
\end{center}
\caption{Hough transform used for detecting geometrical objects. Left: lines detection using MATLAB routines \texttt{houghlines, houghpeaks}. Right: circle detection using MATLAB routine \texttt{imfindcircles}.}
\label{fig:scale_detection}
\end{figure}

\paragraph{Step 4: Measurement scale detection.}  The images in our database divide into two groups: the first is characterised by the presence of linear rulers, whereas the second contains  circular rulers (Figure \ref{fig:im_database}). We thus need to use the Hough transform based Algorithm~\ref{alg:Hough} to detect lines or circles, respectively. The user is then required to tell the software which objects he or she wants to detect. In both cases, in order to avoid false detections (such as ``aligned'' objects erroneously detected as lines, or circle-like objects wrongly considered as circles, see Figure \ref{fig:wrongdetections}), a good candidate for a rough, sensible approximation of the measurement scale is needed as described in Section \ref{subsec:hough}. In order to get this, we proceed as follows: after detecting the head as in Step 1, we use the option \texttt{EquivDiam} of the MATLAB routine \texttt{regionprops} to detect the diameter of the head region (in pixels). We then compare such measurement with pre-collected average measurements of head diameters of male pied flycatchers of a similar population (in $cm$), thus obtaining an initial approximation of the measurement scale. In the case of images containing linear rulers, this will serve as a spacing parameter $s$ for the algorithm. In other words, only lines distant at least $s$ pixels from each other  will be considered. In the case of circular rulers, the same rough approximation will serve similarly as an indication of the range of values in which the Hough transform based MATLAB function \texttt{imfindcircles} will look for circles' radii. For linear ruler images, the algorithm will look only for parallel lines aligned with a prescribed direction. We set this direction as the one perpendicular to the longest line in the image (since the expectation is, that this longest line is the edge of the ruler). Results of this step are shown in Figure \ref{fig:scale_detection}.

\paragraph{Outliers removal for linear rulers.} The scale detection step described above may miss some lines on the ruler. This can be due to an oversmoothing in the denoising step, to high threshold values for edge detection or also to the choice of a large spacing parameter. Furthermore, as we can see from Figures \ref{fig:im_database} and \ref{fig:scale_detection},  we can reasonably assume that the ruler lies on a plane, but its bending can distort some distances between lines. Moreover, few other false line detections can occur (like the number $11$ marked on the ruler main body in Figure \ref{fig:scale_detection}). To exclude these cases, we compute the distance (in pixels) between all the consecutive lines detected and eliminate possible outliers using the standard interquartile range ($IQR$) formula \cite{uptoncook} for outliers' removal. Indicating by $Q_1$ and $Q_3$ the lower quartile and the third quartile, an outlier is every element not contained in the interval $[Q_1-1.5*(Q_3-Q_1), Q_3+1.5*(Q_3-Q_1)]$. Finally, we compute the empirical mean, variance and standard deviation (SD) of the values within this range, thus getting a final indication of the scale of the ruler together with an indicator of the precision of the method. 


\paragraph{Step 5: Measurement.} Once the measurement scale has been detected, it is easy to get all the required measurements. We are interested in the perimeter, the area of the blaze and also in the height and width of the whole blaze component. For linear rulers, due to the error committed in the scale detection step, these values present some uncertainty and variability (see above). In Table \ref{precision} we show the results of numerical tests on a sample of $30$ images with linear rulers. For every image in the sample we compute the standard deviation (SD) error and report in the table the minimum, maximum, and average SD error over the single ones compute, together with the relative standard deviation (RSD) which gives a percentage indication of the error committed.
\begin{equation*}
RSD := \frac{\sigma}{\bar{X}}\cdot 100,
\end{equation*}
where $\sigma$ is the sample SD and $\bar{X}$ is the sample mean of measurements. We observe a minimum and maximum SD of $4.00$ and $10.67$ pixels, respectively, which, compared to the dimension of the original image ($3648\times 2736$ pixels) suggests a reasonable precision. This is confirmed by the average SD value over the sample which is found to be $6.81$ pixels. In percentage, the average error over the sample is $11.99\%$. For circular rulers, we observed in all our experiments that an initial approximation of the range of values for the circle radius (see Step 4 above) results in a robust and typically outlier-free detection of the circular ruler and consequently in an accurate measurement of its radius; the only possible cause of variability and error is its bending.

Uncertainty in the measurements of lengths and areas is calculated with standard formulas in propagation of errors.

\begin{table*}[!h]
\begin{center}
\begin{tabular}{ | c | c | c | c | c | c |}
\hline
 \textbf{SD min} & \textbf{SD max} & \textbf{mean SD}  & \textbf{RSD min} & \textbf{RSD max} & \textbf{mean RSD} \\ \hline
$4.01$ pixels & $10.67$ pixels & $6.81$ pixels & $6.59~ \%$ & $17.36~ \%$ & $11.99~ \%$ \\ \hline
\end{tabular}
\end{center}
\caption{Precision of the measurement scale detection for linear rulers on a sample of $30$ images. The minimum, maximum and average standard deviation (SD) error together with the corresponding relative standard deviation (RSD) errors are reported.}
\label{precision}
\end{table*}

Despite these variabilities, our method is a flexible and semi-supervised approach for this type of problem. Further tests on the whole set of images and improvements on its accuracy are a matter of future research. The analysis of the resulting data measurements for the particular problem of flycatchers' blaze segmentation will be the topic of the forthcoming paper \cite{ourpreprint}. 

We compare in Table \ref{table:measurement} between the use of our combined approach and the use of the manual \emph{line tool} of the IMAGEJ software for the measurement of the blaze area. Namely, we measured in Figure \ref{database:imb} and in Figure \ref{database:imc} the ruler scale by means of the IMAGEJ line tool by considering, for each image, two different $3$ $cm$-sections of the ruler; we then measured manually the number of pixels contained in each, divided each measurement by $30$ and averaged the two results to obtain an estimate of the ruler scale (i.e. the number of pixels crossed by a 1 mm horizontal or vertical line segment). We then measured the area of the blaze after segmenting it by means of the `magic-wand' \cite{morenovelando} IMAGEJ tool and trapezium fitting \cite{pottimontalvo} (see Figure \ref{fig:previous_methods_blaze}). The results are reported  in Table \ref{table:measurement}  both as number of image pixels inside the blaze and in $mm^2$, where this second value has been calculated using the measurement scale detected as described above. We then repeated such measurements using our fully automated Hough transform method for ruler scale detection, reporting as above the measurements of the blaze area computed both as number of image pixels and in $mm^2$. We observe a good level of accuracy of our combined method (see also Table \ref{precision}) with respect to the `magic-wand' manual approach of Moreno \cite{morenovelando}, while, unsurprisingly, the blaze measurements obtained by pure trapezium fitting as proposed by Potti and Montalvo in \cite{pottimontalvo} tend to overestimate the area of the blaze.

\begin{table*}[!h]
\begin{center}
 \begin{tabular}{*{9}{|c}|} 
  \cline{2-9}
  \multicolumn{1}{c}{} & \multicolumn{2}{| c |}{Scale (\emph{\# pixels = 1mm})} & \multicolumn{3}{ | c |}{Blaze area (pixel count)} & \multicolumn{3}{ | c |}{Blaze area ($mm^2$)} \\ 
\cline{2-9}
  \multicolumn{1}{c |}{}    & \textbf{Manual} & \textbf{HT (Ours)} & \textbf{MW} & \textbf{Trap.} & \textbf{GL (Ours)} & \textbf{MW} & \textbf{Trap.} & \textbf{GL (Ours)} \\ \hline
    Figure \ref{database:imb} & $70.2504$   & $72.551$   & $85026$ & $117415$ & $84831$ & $17.2288$ & $23.7917$  & $16.1164$  \\ \hline
   Figure \ref{database:imc} & $71.863$   & $71.8367$  & $101730$ & $146751$ & $121360$ & $19.6980$ & $28.4165$ & $23.517$  \\ \hline
  \end{tabular}
  \end{center}
  \caption{Comparison between ruler scale detection by using manual IMAGEJ line tool and our Hough Transform (HT) method with corresponding measurements of the segmented blaze area obtained by using IMAGEJ `magic-wand' (MW) tool \cite{morenovelando}, trapezium fitting (Trap.) \cite{pottimontalvo} (see also Figure \ref{fig:previous_methods_blaze}) and the graph Ginzburg-Landau (GL) minimisation.}
    \label{table:measurement} 
\end{table*}

\subsection{Moles monitoring for melanoma diagnosis and staging}

In this section we focus on another application of the scale detection Algorithm \ref{alg:Hough} in the context of melanoma (skin cancer) monitoring, see Figure \ref{fig:skinmoles}.  Early signs of melanoma are sudden changes in existing moles and are encoded in the mnemonic ABCD rule. They are \textbf{A}symmetry, irregular \textbf{B}orders, variegated \textbf{C}olour and \textbf{D}iameter \footnote{ \emph{Prevention: ABCD's of Melanoma}. American Melanoma Foundation, \url{http://www.melanomafoundation.org/prevention/abcd.htm}.}.
In the following we focus on the D sign.

Due to their dimensions and their irregular shapes, moles are often very hard to measure. Typically, a common dermatological practice consists in positioning a ruler under the mole and then taking a picture with a professional camera. Sudden changes in the evolution of the mole are then observed by comparison between different pictures taken over time. Hence, their quantitative measurement may be an indication of a malignant evolution   

\begin{figure*}[!h]
\begin{subfigure}{.32\textwidth}
\centering
\includegraphics[height=3cm]{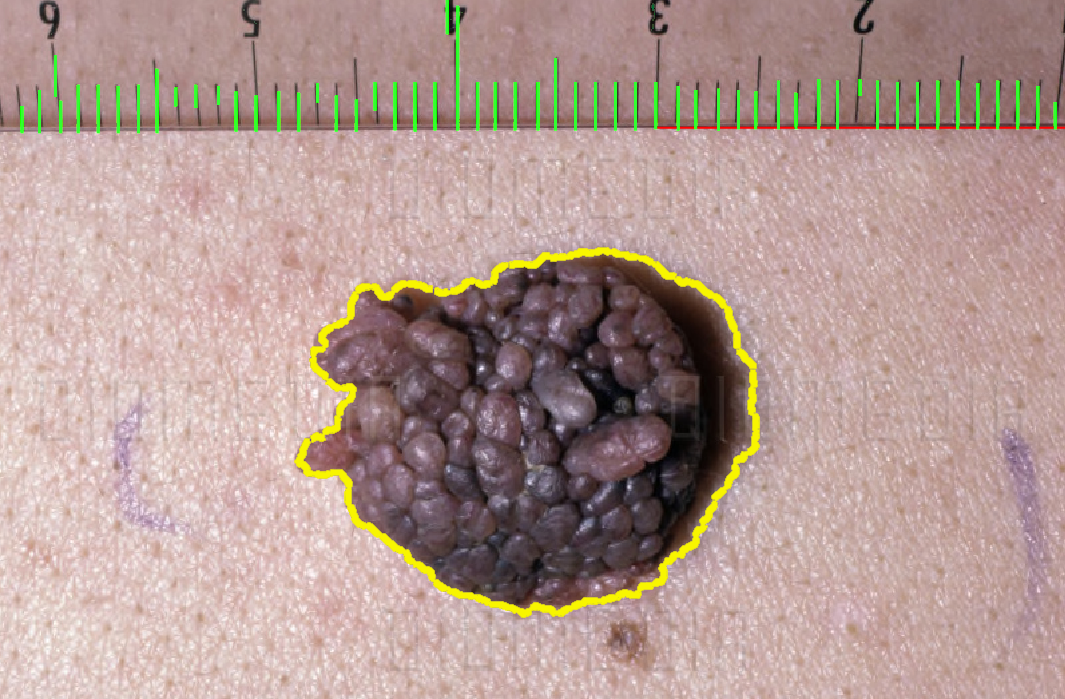}
\caption{}
\label{fig:skinmoles_detection1}
\end{subfigure}
\begin{subfigure}{.32\textwidth}
\centering
\includegraphics[height=3cm]{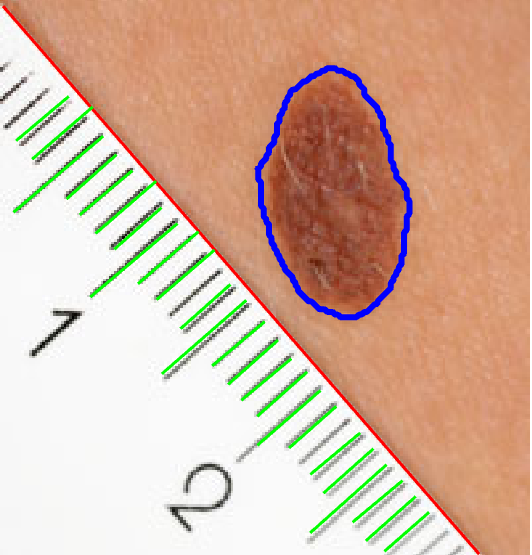}
\caption{}
\label{fig:skinmoles_detection2}
\end{subfigure}
\begin{subfigure}{.32\textwidth}
\centering
\includegraphics[height=3cm]{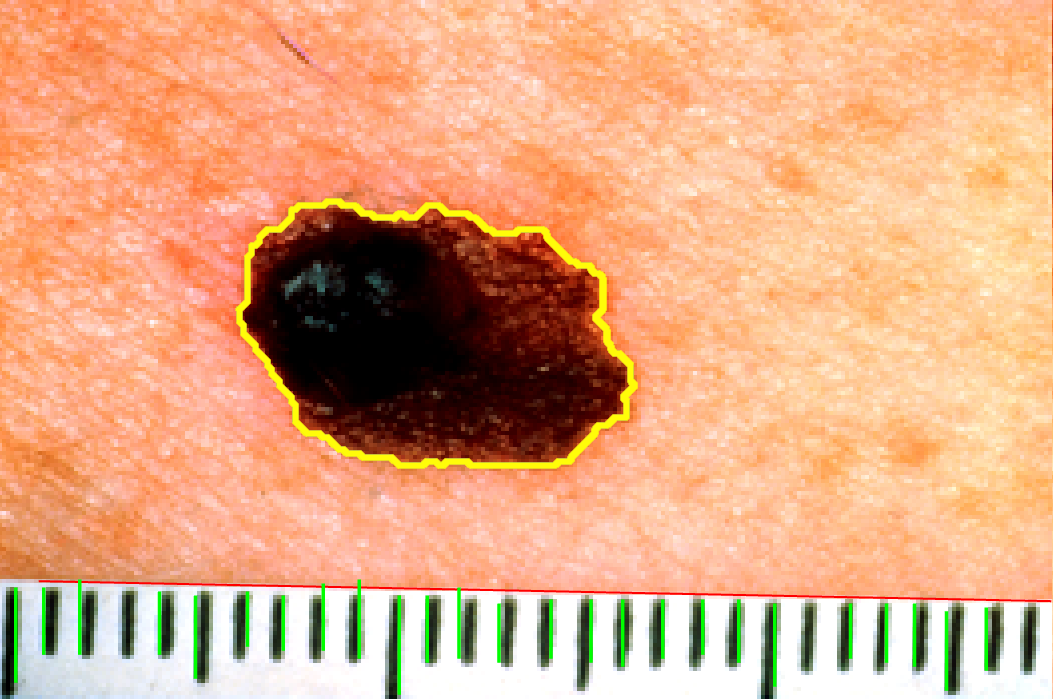}
\caption{}
\label{fig:skinmoles_detection3}
\end{subfigure}
\caption{Moles' detection using GL Algorithm \ref{alg:GLmin_nystrom} (a), the Chan-Vese model \cite{chanvese} ((b),(c)), and measurement scale detection by Hough transform (Algorithm \ref{alg:Hough}).}
\label{fig:skinmoles_detection}
\end{figure*}

In the following examples reported in Figure \ref{fig:skinmoles_detection}, we use the graph segmentation approach described in algorithm \ref{alg:GLmin_nystrom} where texture-characteristic regions are present (see Figure \ref{fig:skinmoles_detection1}) and  the Chan-Vese model \cite{chanvese} for images characterised by homogeneity of the mole and skin regions and the regularity of mole boundaries (Figures \eqref{fig:skinmoles_detection2}-\eqref{fig:skinmoles_detection3}). For the numerical implementation, we use the freely available online IPOL Chan-Vese segmentation code \cite{getreuer}. Let us point out here  that previous works using variational models for accurate melanoma segmentation already exist in literature, see \cite{cavalcanti,abbas}, but in those no measurement technique is considered.

\subsection{Other applications: animal tracks and archeological finds' measurement}

We conclude this section presenting some other applications for the combined segmentation and scale detection models presented above.

\smallskip

\begin{figure*}[!h]
\begin{subfigure}{.5\textwidth}
\begin{center}
\includegraphics[height=3cm]{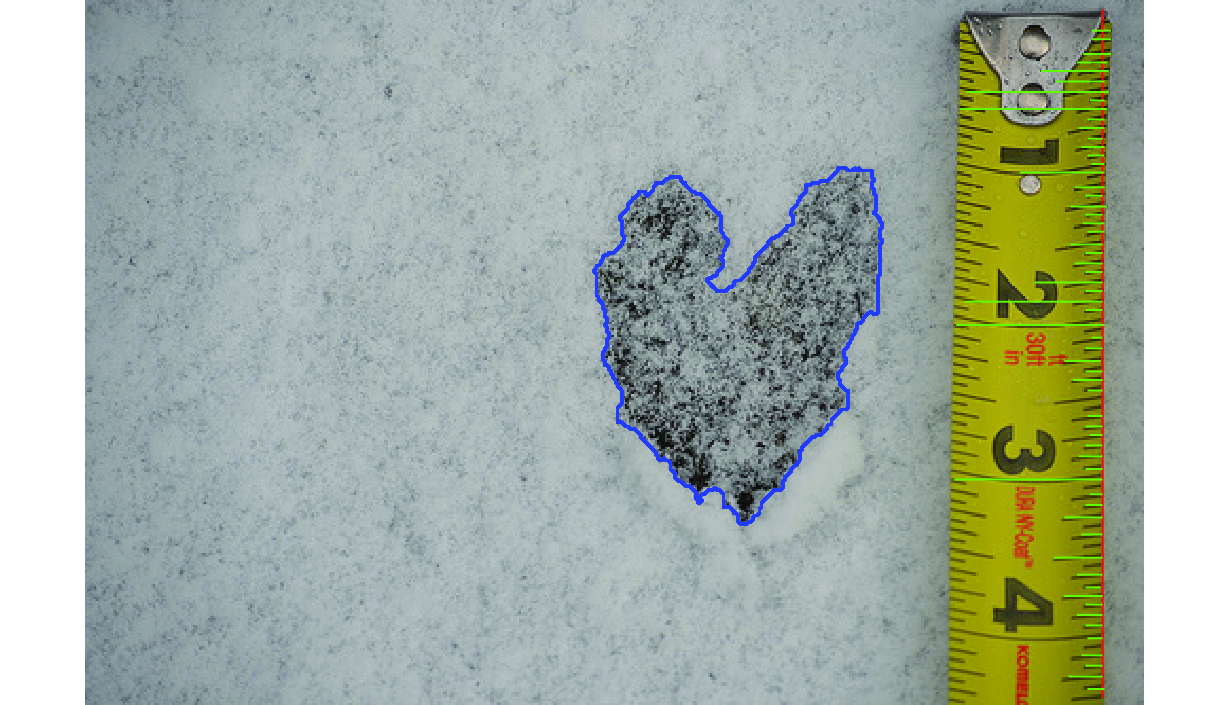}
\end{center}
\caption{White-tailed deer tracks measurement}
\label{fig:footprint}
\end{subfigure}
\hspace{-1cm}
\begin{subfigure}{.5\textwidth}
\begin{center}
\includegraphics[height=3cm]{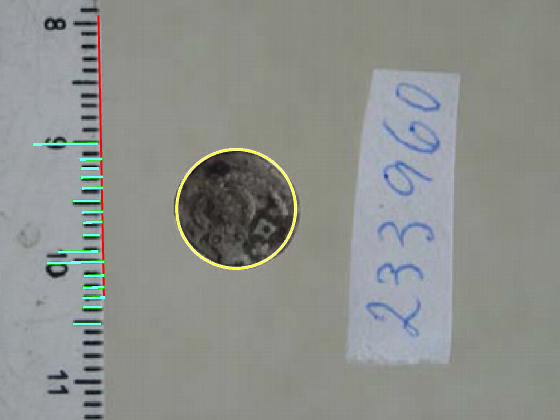}
\end{center}
\caption{Coin measurement, image taken from \cite{herrmann}}
\label{fig:coin}
\end{subfigure}
\caption{The measurement scale has been detected only in a portion of the figure for the sake of reading clarity.}
\label{fig:otherapp}
\end{figure*}

The first application is the identification and classification of animals living in a given area through their soil, snow and mud footprints. Their quantitative measurement is also interesting in the study of the age and size a of a particular animal species. As in the problems above, such measurement very often reduces to a very inaccurate measurement performed with a ruler placed next to the footprint image. In Figure \ref{fig:footprint}\footnote{Image from \url{http://mamajoules.blogspot.co.uk/2015/01/a-naturalists-thoughts-on-animal-tracks.html}.} our combined method is applied for the measurement of a white-tailed deer footprint.

As a final application, we focus on archaeology. In many archaeological finds, objects need to be measured for comparisons and historical studies \cite{herrmann}. Figure \ref{fig:coin} shows the application of our method to coin measurements. Due to its circular shape, for this image a combined Hough transform method for circle and line detection has been used. The example image is taken from \cite{herrmann} where the authors propose a gradient threshold based method combined with a Fourier transform approach. Despite being quite efficient for the particular applications considered, such approach relies in practice on the good experimental setting in which the image is taken: almost noise-free images and very regular objects with sharp boundaries (mainly coins) and homogeneous backgrounds are considered. Furthermore, results are reported only for rulers with vertical orientation and no bending. 


\section{Conclusions}

In this paper we consider image segmentation applications involving measurement of a region's size, which has applications in several disciplines. For example, zoologists may be interested in quantitative measurements of some parts of the body of an animal, such as distinctive regions characterised by specific colours and texture, or in animal tracks to differentiate between individuals in the species. In medical applications, quantifying an evolving, possibly malignant, mass (like, for instance, skin melanoma) is crucial for an early diagnosis and treatment. In archaeology, finds need to be measured and classified. In all these applications, often a common measurement tool is juxtaposed to the region of interest and its measurement is simply read directly from the image. This practice is typically inaccurate and imprecise, due to the conditions in which pictures are taken. There may be noise corrupting the image, the object to be measured may be hard to distinguish, and the measurement tool can be misplaced and far from the object to measure. Moreover, the scale of the image depends on the image itself due to the varying distance from the camera of the ruler and objects to measure.

The method presented (based on \cite{BerFlen}) consists of a semi-supervised approach which, by training the algorithm with some examples provided by the user, extracts relevant features from the training image (such as RGB intensities, texture) and uses them to detect similar regions in the unknown image. Mathematically, this translates into the minimisation of the discrete Ginzburg-Landau functional defined on graphs. To overcome the computational issues due to the size of the data, Nystr{\"o}m matrix completion techniques are used and for the design of an efficient numerical scheme, convex splitting is applied. The measurement scale detection task is performed  by using the Hough transform, a geometrical transformation which is capable of detecting objects with \emph{a priori} known geometrical shapes (like lines on a ruler or circles with fixed diameter). Once the measurement scale is detected, all the measurements are converted into a unit of measure which is not image-dependent. 

Our method represents a systematic and reliable combination of segmentation approaches applied to several real-world image quantification tasks. The use of dictionaries, moreover, allows for flexibility as, whenever needed, the training database can be updated. With respect to recent developments \cite{DataMining} in the fields of data mining for the analysis of big data, where predictions are often performed using training sets and clustering, our approach represents an interesting alternative to standard machine learning (such as $k$-means) algorithms.

\small{\paragraph{Acknowledgements} Many thanks to Colm Caulfield who has introduced HMR and the bird segmentation problem to us mathematicians and to Andrea Bertozzi for her very useful comments on the manuscript. LC acknowledges support from the UK Engineering and Physical Sciences Research Council (EPSRC) grant EP/H023348/1 for the University of Cambridge Centre for Doctoral Training, the Cambridge Centre for Analysis (CCA). CBS acknowledges support from EPSRC grants Nr. EP/J009539/1 and EP/M00483X/1. Moreover, this project has been supported by King Abdullah University of Science and Technology (KAUST) Award No. KUK-I1-007-43. HMR is currently supported by an Institute Research Fellowship at the Institute of Zoology, by the Department of Zoology at the University of Cambridge, and Churchill College, Cambridge. }

\begin{appendices}

\section{The Nystr\"{o}m extension} \label{nystrom:app}

With respect to the eigenvalue problem formulation \eqref{nystrom2} and \eqref{nystrom3}, we revise in this section the Nystr\"{o}m extension \cite{nystrom} in a matrix form.

Let us define first the sub matrices $W_{XX}\in\R^L\times\R^L$ and $W_{XY}\in\R^{L}\times\R^{S-L}$ as
 \begin{align}  \label{matrixWxxWxy}
& W_{XX} =
 \begin{pmatrix}
  w(x_1,x_1) &  \cdots & w(x_1,x_L) \\
  \vdots  &  \ddots & \vdots  \\
  w(x_L,x_1)   & \cdots & w(x_L,x_L)
 \end{pmatrix}, \\
& W_{XY} =
 \begin{pmatrix}
  w(x_1,y_1) &  \cdots & w(x_1,y_{S-L}) \\
  \vdots  &  \ddots & \vdots  \\
  w(x_L,y_1)   & \cdots & w(x_L,y_{S-L})
 \end{pmatrix}. \notag
 \end{align}
Analogous definitions hold for $W_{YY}$ and $W_{YX}$. Each of these matrices represents the sub matrix having as elements the weights between the points in $X$, $Y$ or between the two sets. With this notation, the whole matrix $W\in\R^S\times\R^S$ can be written in block-form as
 \begin{equation*}
 W=
 \begin{pmatrix}
  W_{XX} & W_{XY} \\
  W_{YX}  & W_{YY}
 \end{pmatrix},\qquad W_{YX}=W_{XY}^T.
 \end{equation*}
 Similarly, vectors $v\in\R^S$ can be written as $v=(v_X^T\  v_Y^T)^T$. We focus on the spectral decomposition of the first block of $W$, that is $W_{XX}$. Since this matrix is symmetric, calling $\Theta_X$ the matrix $\Theta_X=diag(\theta_1,\ldots,\theta_L)$ containining the eigenvalues of $W_{XX}$, then by the spectral theorem $W_{XX}=V_X \Theta_X V_X^T$ (compare with \eqref{nystrom3}), with $V_X$ be the orthogonal matrix having as columns the eigenvectors of $W_{XX}$. Writing \eqref{nystrom2} for $y\in Y$, in operator form, we obtain $V_Y$ as
 \begin{equation*}
 V_Y=W_{YX}V_X\Theta_X^{-1}.
 \end{equation*}
 The approximated eigenvectors of $W$ can then be written in matrix form as
 \begin{equation}   \label{approx_eigenvect_W}
 V= \begin{pmatrix}
  V_X \\
 W_{YX}V_X \Theta_X^{-1}
 \end{pmatrix}.
 \end{equation}
Let us observe that
\begin{align}   \label{nystrom_approximation}
 &  V\Theta_X V^T=\begin{pmatrix}
  V_X \\
  W_{YX}V_X \Theta_X^{-1}
 \end{pmatrix} \Theta_X\ [V_X^T\quad (W_{YX}V_X\Theta_X^{-1})^T ] \notag\\ 
 & =\begin{pmatrix}
  V_X\Theta_X V_X^T & W_{XY} \\
  W_{YX} & W_{YX}W_{XX}^{-1}W_{XY}
 \end{pmatrix} = \begin{pmatrix}
  W_{XX} & W_{XY} \\
  W_{YX} & W_{YX}W_{XX}^{-1}W_{XY}
 \end{pmatrix}\approx W. 
\end{align}
Therefore, Nystr\"om extension can be interpreted as the approximation $W\approx V\Theta_X V^T$, under the approximation of $W_{YY}$ given by
$
W_{YY}\approx W_{YX}W_{XX}^{-1}W_{XY}.
$
The quality of the approximation of the full $W$ is quantified by the norm of the Schur complement $\| W_{YY}- W_{YX}W_{XX}^{-1}W_{XY}\|$, see \cite{fowlkes} 

Recalling the definition of the symmetric graph Laplacian $L_s$ given by \eqref{symmgraphlapl} and the relation between the spectral decomposition of $W$ and the one of $W$ in \eqref{eigenvectorsLs}, we observe that a normalisation step now needs to be computed for obtaining the spectral decomposition of $L_s$. Defining $\textbf{1}_L$ as the $L$-dimensional vector consisting of ones and $\textbf{1}_{S- L}$ analogously, we use \eqref{nystrom_approximation} and start computing the nonnegative vector $d=(d_X^T d_Y^T)^T$ by
\begin{align}  \label{nystrom_normalise1}
& \begin{pmatrix}
d_X \\
d_Y
\end{pmatrix} =
\begin{pmatrix}
  W_{XX} & W_{XY} \\
  W_{YX} & W_{YX}W_{XX}^{-1}W_{XY}
 \end{pmatrix}~
 \begin{pmatrix}
\textbf{1}_L \\
\textbf{1}_{S-L}
\end{pmatrix}=
\begin{pmatrix}
  W_{XX}\textbf{1}_L+ W_{XY}\textbf{1}_{S-L}  \\
  W_{YX}\textbf{1}_L + W_{YX}W_{XX}^{-1}W_{XY}\textbf{1}_{S-L}
 \end{pmatrix}. 
\end{align}
Therefore, the matrices $W_{XX}$ and $W_{XY}$ can be normalised simply by considering:
\begin{align}   \label{nystrom_normalise2}
& \hat{W}_{XX}=W_{XX}./(\sqrt{d_X}\otimes\sqrt{d_X}^T), \\ 
& \hat{W}_{XY}=W_{XY}./(\sqrt{d_Y}\otimes\sqrt{d_Y}^T),  \notag
 \end{align}
where the division is intended element-wise and $\otimes$ is the standard vector tensor product.

A further step of normalisation is now needed since the approximated eigenvectors of $W$, i.e. the columns of the matrix $V$ in \eqref{approx_eigenvect_W} may not be orthogonal. Such normalisation may be obtained by using auxiliary unitary matrices. We refer the reader to \cite[Section 3.2]{BerFlen} for more details on this. 

Once these additional normalisation steps are completed, we then get a spectral decomposition of $W$ in terms of its eigenvalues $\hat{\lambda}_i$ and the corresponding normalised eigenvectors $v_i,~i=1,\ldots,S$. Therefore, recalling \eqref{eigenvectorsLs}, the spectral decomposition of $L_s$ is given in terms of the eigenvalue $1-\hat{\lambda}_i$ and eigenvectors $v_i$.

\red{\section{The MBO scheme for image segmentation} \label{appendix:MBO}
As previously commented in Section \ref{subsec:GL}, by taking the $L^2$ gradient descent of the Ginzburg-Landau functional defined in \eqref{def:GL_func}, one gets the well-known Allen-Cahn equation \cite{alcahn}:
\begin{equation}   \label{eq:al_cahn}
u_t= \varepsilon\Delta u - \frac{1}{\varepsilon}W'(u),
\end{equation}
which has often been studied for the modelling of several phase transition and separation problems and for the study of mean curvature flow (see, e.g., \cite{brokate}). In the limit $\epsilon \to 0$ solutions consist of two phases corresponding to the wells of $W$. In \cite{Rubinstein} it is shown that, for rescaled solutions of equation \eqref{eq:al_cahn}, the interface between these phases evolves according to mean curvature flow. In \cite{MBO1992}, Merriman, Bence and Osher propose an alternative approach (later named MBO scheme) which, by using threshold dynamics, approximates the mean curvature flow of the interface at discrete times. As proved rigorously in \cite{barles}, for small values of the interface parameter $\varepsilon$, the MBO scheme can then be used to solve equation \eqref{eq:al_cahn} numerically. }

\red{
In \cite{MerkujevKostic2013}, the authors propose a variant of the MBO scheme as an alternative way to (approximately) minimise the graph GL functional with fidelity term, \eqref{def:GL_discr}.
Recalling the graph framework introduced in Section \ref{subsec:graph}, the MBO segmentation starts from an initialisation $U_1$ given by \eqref{initialisation:GL} and computes, for every $n\geq 1$ the new iterate $U_{n+1}$ from $U_n$ by applying sequentially the two following steps:
\begin{itemize}
\item \textbf{Step 1} (diffusion with forcing term): Starting from $U^1_n=U_n$, solve for every $1\leq k \leq K$ the discretised heat diffusion equation with fidelity term
\begin{equation}  \label{MBO:step1}
\frac{U^{k+1}_n-U^{k}_n}{\tau} = -L_s~U^{k+1}_n-\chi(x)(U^{k+1}_n-U_0),
\end{equation}
where $\tau:=\frac{\Delta t}{K}$ is the heat diffusion time and $K$ is the number of diffusion steps. Practically, $\tau$ has to be chosen small enough to approximate the motion by mean curvature and large enough to avoid \textit{freezing} or \textit{pinning}, which occurs when the diffusion time is so short that not enough mass diffuses along the edges of the network and the thresholding operation described in the following Step 2 leaves $U_n$ unchanged.
\item \textbf{Step 2} (thresholding): For every point $x$ set $U_{n+1}$ as:
\begin{equation*} 
U_{n+1}(x)=\begin{cases}
1,\quad&\text{if }U^K_n(x)\geq 0, \\
-1,\quad&\text{if }U^K_n(x)< 0.
\end{cases}
\end{equation*}
\end{itemize}
Numerically, \eqref{MBO:step1} is solved at each diffusion time step $k\tau, ~k\geq 1$ by considering the spectral decomposition of $U_n^k$ with respect of the eigenvectors of the operator $L_s$, similarly as in \eqref{spectrdecun}, and using classical Fourier transform methods to compute the new iterate $U_n^{k+1}$.
}

\section{The Hough transform} \label{sec:hough}

The general idea behind the use of the Hough transform 
\cite{hough,dudahart} is to map the ambient space to an auxiliary space, called the parameter space (as it is related to the parametric representation of the geometrical objects we are interested in). There, objects with specified shapes are easily recognisable as peaks of specific functions. Let us clarify these concepts with two examples.

\begin{figure}[h!]
\begin{center}
\begin{subfigure}{.2\textwidth}
\centering
\includegraphics[height=3.8cm]{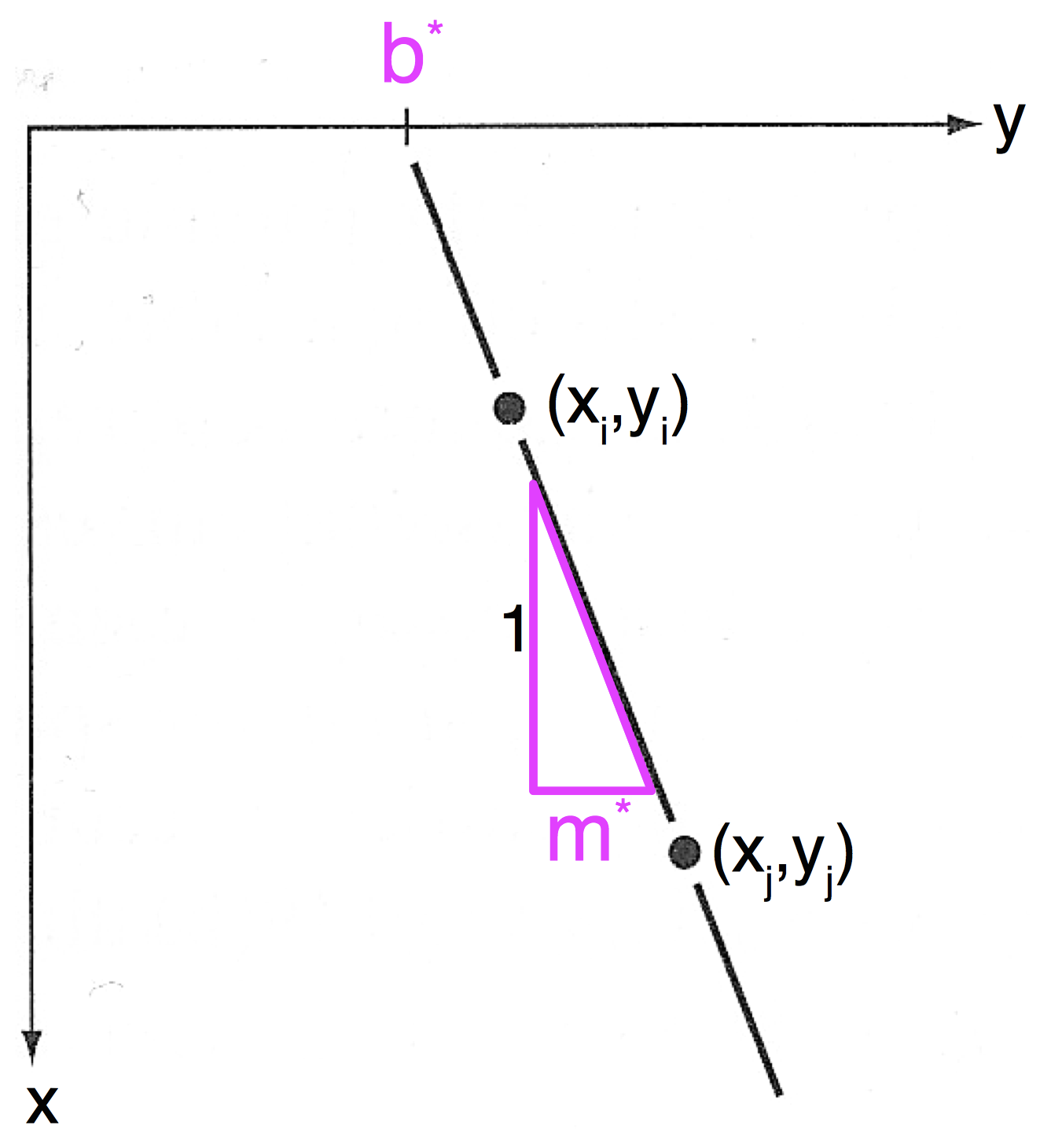}
\caption{$x$-$y$-plane}
\label{fig:houghxy}
\end{subfigure}
\begin{subfigure}{.2\textwidth}
\centering
\includegraphics[height=3.8cm]{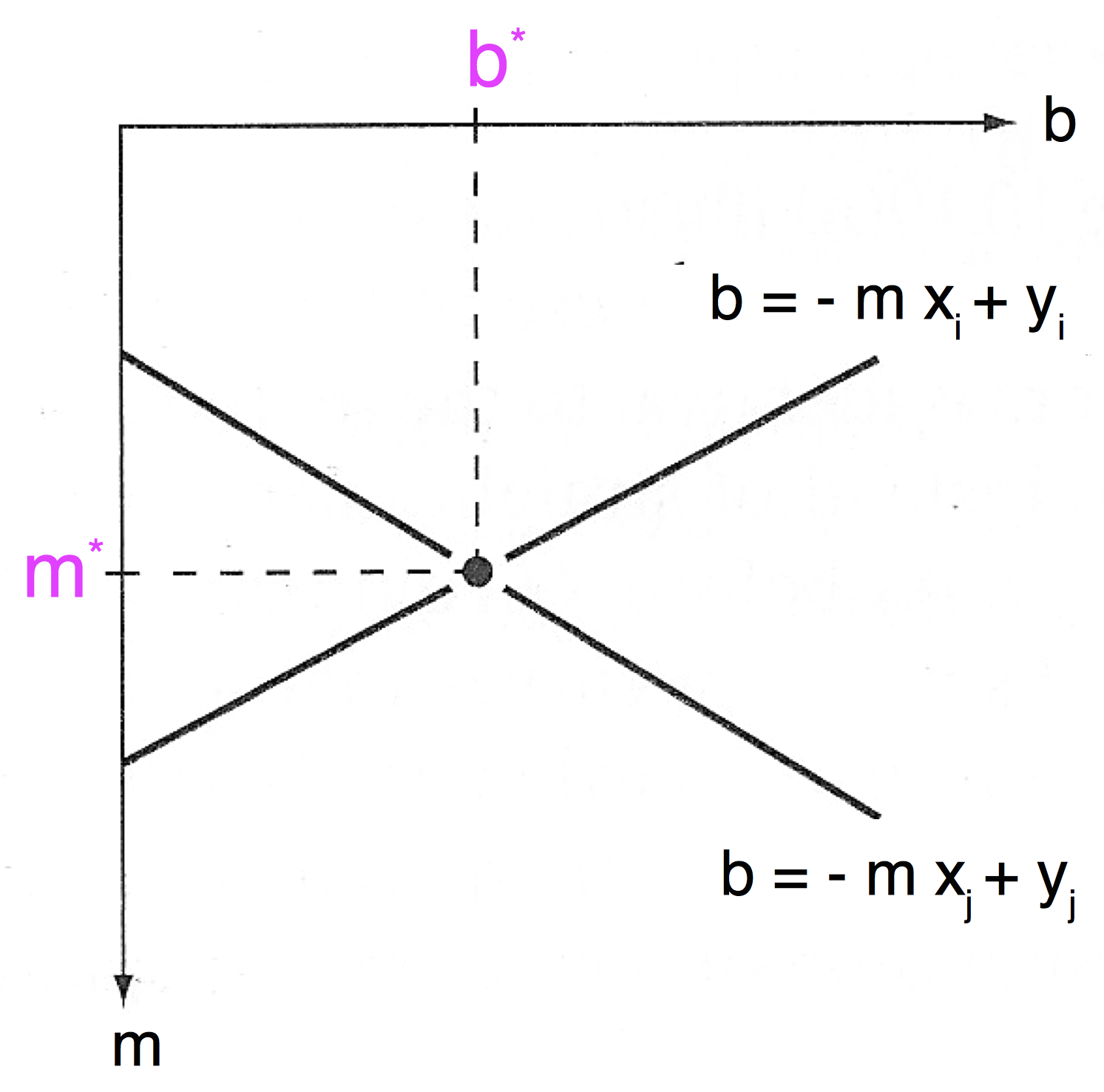}
\caption{$m$-$b$-plane}
\label{fig:houghmq}
\end{subfigure}
\caption{Slope-intercept form, \eqref{slope_intercept}. Images edited from \cite{Gonzalez}.}
\label{fig:HoughSlopeIntercept}
\end{center}
\end{figure}

 \paragraph{Detecting line segments.} We start from the typical \emph{slope-intercept form} of a line:
 \begin{equation}\label{slope_intercept}
 y=mx+b, \quad m,b,x,y\in\R.
 \end{equation}
Traditionally, the equation above is considered as a function of points with coordinates $(x,y)$ satisfying equation \eqref{slope_intercept} for \emph{fixed} values of $m$ and $b$. In other words, these values identify a specific straight line in the $x$-$y$ plane, cf. Figure \ref{fig:houghxy}. Rewriting \eqref{slope_intercept} as $b=y-mx$ and keeping fixed the coordinates $(x,y)$ we obtain a new equation of a straight line in the $m$-$b$ plane, cf. Figure \ref{fig:houghmq}, depicting the parameter space. If lines in the $m$-$b$ parameter space intersect, their sign-changed slopes (given by their $x$ values) and $m$-intercepts (their $y$ values) correspond to points lying on the same line in the $x$-$y$ plane. The $(m,b)$ coordinates of the intersection point in parameter space specify the slope and $x$-intercept respectively of that line in the $x$-$y$ plane.

Hence, if we are given a black and white image in the $x$-$y$ plane, and for all coordinates $(x,y)$ of black locations in the image, we draw the corresponding lines in the $m$-$b$ plane, intersection points of those lines will tell us which $(x,y)$ locations in the image lie on the same line. Of course any two points lie on a line, thus we are specifically interested in intersection points in the $m$-$b$ plane in which many different lines intersect, indicating the presence of an actual black line segment in the original image.

\begin{figure*}[h!]
\begin{subfigure}{.32\textwidth}
\centering
\includegraphics[height=3.5cm]{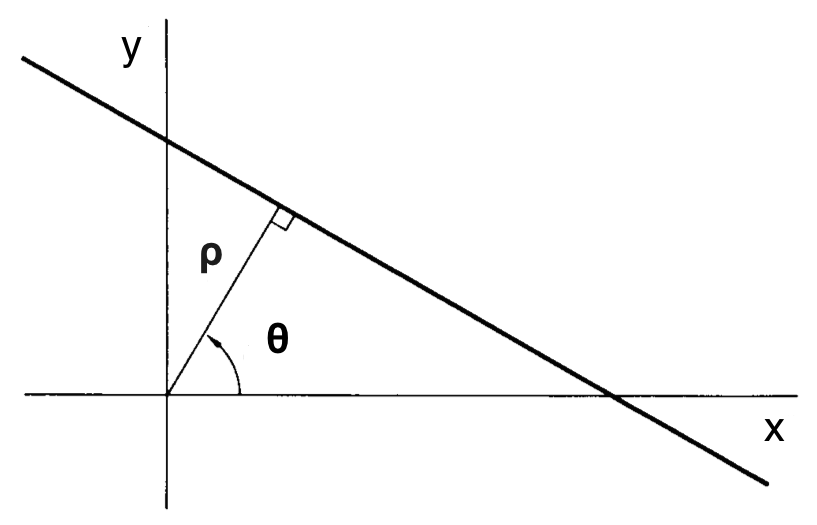}
\caption{Normal parametrisation \eqref{normal_param}}
\label{fig:normalparam}
\end{subfigure}
\begin{subfigure}{.32\textwidth}
\centering
\includegraphics[height=3.5cm]{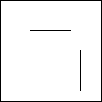}
\caption{Binary image with two lines}
\label{fig:binim}
\end{subfigure}
\begin{subfigure}{.32\textwidth}
\centering
\includegraphics[height=3.5cm]{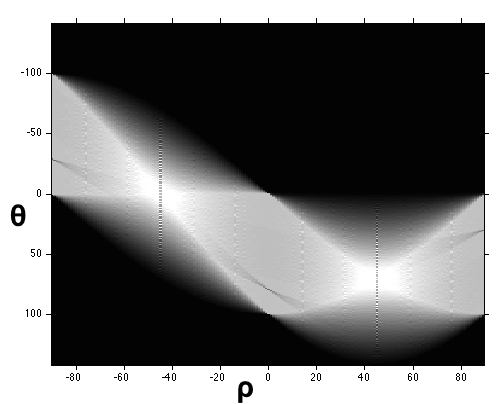}
\caption{Parameter space}
\label{fig:binimparsp}
\end{subfigure}
\caption{Normal form: image and parameter spaces.}
\label{fig:HoughNormal}
\end{figure*}

Drawbacks of this parametrisation are the need for an unbounded parameter space to describe near vertical lines and the impossibility to describe a vertical line. One alternative is the \emph{normal parametrisation} which views a straight line in $x$-$y$ space as the tangent line to a circle with radius $\rho$, touching the circle at angular coordinate $\theta$, as illustrated as in Figure \ref{fig:normalparam}, \cite{dudahart}. In $\rho$-$\theta$ parameter space this leads to
\begin{equation}   \label{normal_param}
\rho=x\cos\theta+y\sin\theta, \quad \theta\in [0,\pi].
\end{equation}

The objects in the parameter space are now sinusoidal curves, but again intersection points identify parameters for the points lying on the same straight line in the $x$-$y$ plane. Figures \ref{fig:binim} and  \ref{fig:binimparsp} show a binary image with two black straight lines and the corresponding parameter space. The bright spots in the parameter space indicate a large number of intersections, thus identifying the two lines in the original image. 


\begin{figure*}[h!]
\begin{subfigure}{.45\textwidth}
\centering
\includegraphics[height=4cm]{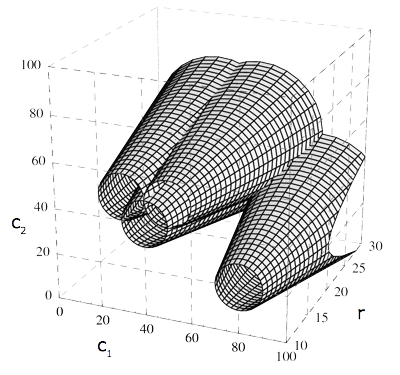}
\caption{Cones in 3D parameter space}
\label{fig:cones}
\end{subfigure}
\begin{subfigure}{.45\textwidth}
\centering
\includegraphics[height=4cm]{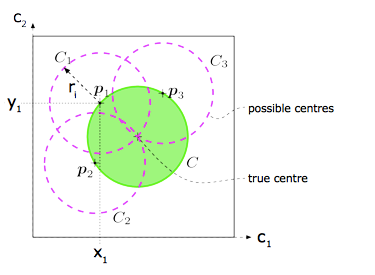}
\caption{Slice of 3D parameter space}
\label{fig:accarray}
\end{subfigure}
\caption{Circular Hough transform. Images edited from \cite{burgerburge}.}
\label{fig:CircularHough}
\end{figure*}

\paragraph{Detecting circles.} Analogously to what we did above, when looking for circular structures in a given image, we consider, for $(x,y)\in\R^2$ the parametric representation of a circle,
\begin{equation}
\label{eq:Circle}
r^2=(x-c_1)^2+(y-c_2)^2,
\end{equation}
where $r>0$ is the radius of the circle and $(c_1, c_2)\in \R^2$ are the coordinates of its centre. Every point $(x,y)$ lying on the circle, satisfies equation \eqref{eq:Circle} for fixed $r, c_1$ and $c_2$. As before, we now consider equation \eqref{eq:Circle} in the three-dimensional parameter space $c_1-c_2-r$ for fixed $x$ and $y$. Here, the objects of interest are cone-shaped surfaces, as shown in Figure \ref{fig:cones}. Their intersection identifies the desired values of $r, c_1$ and $c_2$ in equation \eqref{eq:Circle}, see Figure \ref{fig:accarray}.

\end{appendices}



\end{document}